\documentstyle[11pt,epsfig,epsf]{article}
\input{amssymb.sty}

\begin{document}
\title{\bf Spectra of Upper-triangular Operator Matrix\footnote{This work is
supported by the NSF of China (Grant Nos. 10771034, 10771191 and
10471124) and the NSF of Fujian Province of China (Grant Nos.
Z0511019, S0650009).}}
\author{{Shifang Zhang$^1$, \, Huaijie Zhong$^2$, \, Junde Wu$^1$\footnote{Corresponding
author: E-mail: wjd@zju.edu.cn}}
\\\\$^1$\small\it Department of Mathematics, Zhejiang University,
Hangzhou 310027, P. R. China
\\\\$^2$\small\it Department of Mathematics, Fujian Normal University, Fuzhou 350007, P. R. China}
\date{} \maketitle
\begin{center}
\begin{minipage}{140mm}
{\bf Abstract} {Let $X$ and $Y$ be Banach spaces, $A\in B(X)$, $B\in
B(Y)$, $C\in B(Y, X)$,
$M_{C}=\left(\begin{array}{cc}A&C\\0&B\\\end{array} \right)$ be the
operator matrix acting on the Banach space $X\oplus Y$. In this
paper, we give out 20 kind spectra structure of $M_C$, decide 18
kind spectra filling-in-hole properties of $M_C$, and present 10
examples to show that some conclusions about the spectra structure
or filling-in-hole properties of $M_C$ are not true.}

\vspace{4mm} {\bf Keywords}: \,{Banach spaces, Upper-triangular
operator matrix, Spectra, filling-in-hole.}
\end{minipage}
\end{center}

\vskip 0.2in
%
%
 \large\section{Introduction and basic concepts}

\par It is well known that if $H$ is a Hilbert space and $T$ is a bounded linear operator defined on $H$ and $H_1$ is
an invariant closed subspace of $T$, then $T$ can be represented for
the form of $$T=\left(
       \begin{array}{cc}
        *&*\\
         0&*\\
       \end{array}
     \right):H_1\oplus H_1^{\perp}\rightarrow H_1\oplus H_1^{\perp},$$
which motivated the interest in $2 \times 2$ upper-triangular
operator matrices (see [2], [3], [6], [8-25], [28-32]). Throughout
this paper, let $X$ and $Y$ be complex infinite dimensional Banach
spaces and $B(X, Y)$ be the set of all bounded linear operators from
$X$ into $Y$, for simplicity, we write $B(X, X)$ as $B(X)$. Let
$X^*$ be the dual space of $X$. If $T\in B(X, Y)$, then $T^*\in
B(Y^*,
  X^*)$ denotes the dual operator of $T$.

\par For $T\in B(X, Y)$, let $R(T)$ and $N(T)$ denote the
range and kernel of $T$, respectively, and denote $\alpha (T)=\dim
N(T)$, $\beta (T)=\dim Y / R(T)$. If $T\in B(X)$, the ascent
$asc(T)$ of $T$ is defined to be the smallest nonnegative integer
$k$ (if it exists) which satisfies that $N(T^{k})=N(T^{k+1})$. If
such $k$ does not exist, then the ascent of $T$ is defined as
infinity. Similarly, the descent $des(T)$ of $T$ is defined as the
smallest nonnegative integer $k$ (if it exists) for which
$R(T^{k})=R(T^{k+1})$ holds. If such $k$ does not exist, then
$des(T)$ is defined as infinity, too. If the ascent and the descent
of $T$ are finite, then they are equal (see [13]). For $T\in B(X)$,
if $R(T)$ is closed and $\alpha (T)<\infty$, then $T$ is said to be
an upper semi-Fredholm operator, if $\beta (T)<\infty$, then $T$ is
said to be a lower semi-Fredholm operator. If $T\in B(X)$ is either
upper or lower semi-Fredholm operator, then $T$  is said to be a
semi-Fredholm operator. For semi-Fredholm operator $T$, its index
ind $(T)$ is defined as ind $(T)=\alpha(T )-\beta(T).$

\vspace{3mm}

Now, we introduce the following important operator classes:

\vspace{3mm}

The sets of all invertible operators, bounded below operators,
surjective operators, left invertible operators,  right invertible
operators on $X$ are defined, respectively, by
\begin{eqnarray*}
& & G(X):=\{T \in  B(X):T \makebox{ is invertible}\}, \\
& & G_+(X):=\{T \in  B(X):T \makebox{ is injective and}\ R(T)\makebox{ is closed}\}, \\
& & G_-(X):=\{T \in  B(X):T \makebox{ is surjective}\},\\
& & G_{l}(X):=\{T \in  B(X): T \makebox{ is left invertible}\}, \\
& &G_r(X):=\{T \in  B(X): T \makebox{ is right invertible}\}.
\end{eqnarray*}

The sets of all Fredholm operators, upper semi-Fredholm operators,
lower semi-Fredholm operators, left semi-Fredholm operators, right
semi-Fredholm operators on $X$ are defined, respectively, by
\begin{eqnarray*}
& & \Phi(X):=\{T \in  B(X):\alpha (T)<\infty \makebox{ and }\beta (T)< \infty \},\\
& & \Phi_+(X):=\{T \in  B(X):\alpha (T)<\infty \makebox{ and } R(T)\makebox{ is closed} \}, \\
& &\Phi_-(X):=\{T \in  B(X): \beta (T)< \infty \}, \\
& & \Phi_{l}(X):=\{T \in  B(X):
 R(T)\makebox{ is a closed and complemented subspace of }X  \makebox{ and }\,\,\alpha (T)<\infty\},\\
& &\Phi_r(X):=\{T \in  B(X):
 N(T)\makebox{ is a closed and complemented subspace of }X  \makebox{ and }\,\,\beta (T)<\infty \}.
\end{eqnarray*}

 The sets of all Weyl operators, upper semi-Weyl operators, lower semi-Weyl operators, left semi-Weyl operators,
right semi-Weyl operators on $X$ are defined, respectively, by
\begin{eqnarray*}
& \Phi_{0}(X):=\{T \in \Phi(X): \makebox {ind} (T)= 0 \},\\
& \Phi_+^-(X):=\{T \in  \Phi_+(X):\makebox {ind} (T)\leq 0 \}, \\
 &\Phi_-^+(X):=\{T \in \Phi_-(X):\makebox {ind} (T)\geq 0\}, \\
 &  \Phi_{lw}(X):=\{T \in \Phi_l(X): \makebox {ind} (T)\leq 0 \},\\
 &\Phi_{rw}(X):=\{T \in  \Phi_r(X):\makebox {ind}  (T)\geq 0 \}.
\end{eqnarray*}

The sets of all Browder operators, upper semi-Browder operators,
lower semi-Browder operators, left semi-Browder operators, right
semi-Browder operators on $X$ are defined, respectively, by
\begin{eqnarray*}
& & \Phi_{b}(X):=\{T \in  \Phi(X): asc(T)= des(T)<\infty \},\\
& & \Phi_{ab}(X):=\{T \in  \Phi_+(X): asc(T)<\infty \}, \\
& & \Phi_{sb}(X):=\{T \in  \Phi_-(X): des(T)<\infty \},\\
& &  \Phi_{lb}(X):=\{T \in  \Phi_l(X): asc(T)<\infty \},\\
& &\Phi_{rb}(X):=\{T \in  \Phi_r(X): des(T)<\infty \}.
\end{eqnarray*}

By the help of above set classes, for $T\in B(X)$, we can define its
corresponding 22 kind spectra, respectively, as following:

\vspace{3mm}

\par \par the spectrum: $\sigma_{}(T)=\{\lambda \in
{\mathbb{C}}:T-\lambda I\not \in G(X) \}$,
\par the approximate point spectrum: $\sigma_{a}(T)=\{\lambda \in {\mathbb{C}}:T-\lambda I\not \in G_-(X)\}$,
\par the defect spectrum: $\sigma_{su}(T)=\{\lambda \in {\mathbb{C}}:T-\lambda I\not \in G_+(X)\}$,
\par the left spectrum: $\sigma_{l}(T)=\{\lambda \in {\mathbb{C}}:T-\lambda I\not \in G_l(X)\}$,
\par the right spectrum: $\sigma_{r}(T)=\{\lambda \in {\mathbb{C}}:T-\lambda I\not \in G_r(X)\}$,
\par the essential spectrum:   $\sigma_{e}(T)=\{\lambda \in {\mathbb{C}}:T-\lambda I\not \in \Phi(X)\}$,
\par the upper semi-Fredholm spectrum: $\sigma_{SF+}(T)=\{\lambda \in {\mathbb{C}}: T-\lambda I \not \in \Phi_+(X)\},$
\par the lower semi-Fredholm spectrum:  $\sigma_{SF-}(T)=\{\lambda \in {\mathbb{C}}: T-\lambda I \not \in \Phi_-(X)\},$
\par the left semi-Fredholm spectrum: $\sigma_{le}(T)=\{\lambda \in {\mathbb{C}}: T-\lambda I \not \in \Phi_l(X)\},$
\par the right semi-Fredholm spectrum: $\sigma_{re}(T)=\{\lambda \in {\mathbb{C}}: T-\lambda I \not \in \Phi_r(X)\},$
\par the Weyl spectrum: $\sigma_{w}(T)=\{\lambda \in {\mathbb{C}}:T-\lambda I \not \in \Phi_{0}(X) \},$
\par the upper semi-Weyl spectrum: $\sigma_{aw}(T)=\{ \lambda \in {\mathbb{C}}: T-\lambda I \not \in \Phi_+^-(X)\},$
\par the lower semi-Weyl spectrum: $\sigma_{sw}(T)=\{ \lambda \in {\mathbb{C}} : T-\lambda I \not \in \Phi_{-}^+(X)\},$
\par the left semi-Weyl spectrum: $\sigma_{lw}(T)=\{ \lambda \in {\mathbb{C}}: T-\lambda I \not \in \Phi_{lw}(X)\},$
\par the right semi-Weyl spectrum: $\sigma_{rw}(T)=\{ \lambda \in {\mathbb{C}} : T-\lambda I \not \in \Phi_{rw}(X)\},$
\par the Browder spectrum: $\sigma_{b}(T)=\{\lambda \in {\mathbb{C}}:T-\lambda I \not \in \Phi_{b}(X)\},$
\par the Browder essential approximate point spectrum: $\sigma_{ab}(T)=\{\lambda \in {\mathbb{C}} : T-\lambda I \not \in \Phi_{ab}(X)\},$
\par the lower semi-Browder spectrum: $\sigma_{sb}(T)=\{\lambda \in {\mathbb{C}} : T-\lambda I \not \in \Phi_{sb}(X)\},$
\par the left  semi-Browder spectrum: $\sigma_{lb}(T)=\{ \lambda \in {\mathbb{C}} : T-\lambda I \not \in \Phi_{lb}(X)\},$
\par the right semi-Browder spectrum: $\sigma_{rb}(T)=\{\lambda \in {\mathbb{C}} : T-\lambda I \not \in \Phi_{rb}(X)\},$
\par the Kato spectrum: $\sigma_{K}(T)=\{\lambda \in {\mathbb{C}}: T-\lambda I \not \in \Phi_+(X)\cup \Phi_-(X)\},$
\par the third Kato spectrum: $\sigma_{K_3}(T)=\{\lambda \in {\mathbb{C}}: T-\lambda I \not \in \Phi_l(X)\cup \Phi_r(X)\}.$

\vspace{3mm}

It is well known that all of these spectra are compact nonempty
subsets of complex plane ${\mathbb{C}}$ and have the following
relationship:

\vspace{3mm}

\par(1) $\sigma_{K}(T)\subseteq
\sigma_{SF+}(T)\subseteq\sigma_{aw}(T)\subseteq\sigma_{ab}(T)\subseteq\sigma_{b}(T),$
\par(2) $ \sigma_{K}(T)\subseteq
\sigma_{SF-}(T)\subseteq\sigma_{sw}(T)\subseteq\sigma_{sb}(T)\subseteq\sigma_{b}(T),$
\par (3) $\sigma_{K_3}(T)\subseteq
\sigma_{le}(T)\subseteq\sigma_{lw}(T)\subseteq\sigma_{lb}(T)\subseteq\sigma_{b}(T),$
\par  (4) $\sigma_{K_3}(T)\subseteq
\sigma_{re}(T)\subseteq\sigma_{rw}(T)\subseteq\sigma_{rb}(T)\subseteq\sigma_{b}(T),$
\par (5) $\partial(\sigma_{b}(T))\subseteq
\partial(\sigma_{w}(T))\subseteq\partial(\sigma_{e}(T))\subseteq\sigma_{K}(T)\subseteq
\sigma_{e}(T)\subseteq\sigma_{w}(T)\subseteq\sigma_{b}(T)\subseteq\sigma_{}(T),$
\par (6) $\partial(\sigma_{}(T))\subseteq \sigma_{a}(T)\cap \sigma_{su}(T)\subseteq
\sigma_{l}(T)\subseteq\sigma_{r}(T) \subseteq\sigma_{}(T).$

\vspace{3mm}

\par For a compact subset $M$ of $\mathbb{C}$, we use
$acc M$, $int M$ and $iso M$,
 respectively, to denote all the
 points of accumulation of $M$, the interior of $M$ and
 all the isolated points of $M$.

\vspace{3mm}

An operator $T \in B(X)$ is said to be Drazin invertible if there
exists an operator $T^D \in B(X)$ such that
$$TT^D=T^D T ,\quad \quad T^DTT^D=T^D ,\quad \quad T^{k+1}T^D=T^{k}$$
for some nonnegative integer $k$ ([13], [31]).

\vspace{3mm}

The operator $T^D$ is said to be a Drazin inverse of $T$. It follows
from [13] that $T^D$ is unique. The smallest $k$ in the previous
definition is called as the Drazin index of $T$ and denoted by
$i(T)$. Now, we can define the Drazin spectrum, the ascent spectrum
and the descent spectrum of $T$, respectively, as following:

 \vspace{3mm}

 $$\sigma_{D}(T)=\{\lambda \in {\mathbb{C}}:T-\lambda I \makebox{ is not Drazin invertible} \},$$
 $$\sigma_{asc}(T)=\{\lambda \in {\mathbb{C}}: asc(T-\lambda I )=\infty
 \},$$
 $$\sigma_{des}(T)=\{\lambda \in {\mathbb{C}}: des(T-\lambda I) =\infty \}.$$

\vspace{3mm}

 \par \noindent The sets $\sigma_{D}(T)$, $\sigma_{asc}(T)$ and $\sigma_{des}(T)$
 are closed but may be empty ([7], [31]).

\vspace{3mm}

Now, we continue to introduce the following operator classes which
were discussed in [1], [3-5] and [26-27]:
\begin{eqnarray*}
& &BF(X)=\{T\in B(X):  T=T_1\oplus T_2, \, \makebox{ where}\, T_1 \,\makebox{is a Fredholm operator and }\, T_2\, \makebox{nilpotent}\},\\
& &BW(X)=\{T\in B(X):  T=T_1\oplus T_2, \, \makebox{ where}\, T_1
\,\makebox{is  a Weyl operator and }\, T_2\, \makebox{nilpotent}\},\\
& &R_4(X)=\{T\in B(X): des(T)<\infty,  \,\, R(T^{des(T)})\makebox{ is closed}\},\\
& &R_9(X)=\{T\in B(X): asc(T)<\infty,  \,\, R(T^{asc(T)+1})\makebox{ is closed}\},\\
& & SF_0(X)=\{T\in \Phi_+(X)\cup \Phi_-(X):\alpha(T)=0\,\, \makebox{or} \,\, \beta(T)=0\},\\
& &M(X)=\{T\in B(X): T \,\, \makebox{is left or right invertible}\},\\
& &D(X)=\{T\in B(X): R(T) \,\, \makebox{is closed and} N(T)\subseteq
\cap_{n=1}^\infty R(T^n)\}.
\end{eqnarray*}

\vspace{3mm}

Their corresponding spectra can be defined, respectively, by

\vspace{3mm}

\par The B-Fredholm  spectrum: $\sigma_{BF}(T)=\{\lambda\in{\mathbb{C}}:  T-\lambda I\not\in BF(X)\},$
\par The B-Weyl spectrum: $\sigma_{BW}(T)=\{\lambda\in{\mathbb{C}}:  T-\lambda I\not\in BW(X)\},$
\par The right Drazin spectrum: $\sigma_{rD}(T)=\{\lambda\in{\mathbb{C}}:  T-\lambda I\not\in R_4(X)\},$
\par The left Drazin spectrum: $\sigma_{lD}(T)=\{\lambda\in{\mathbb{C}}:  T-\lambda I\not\in R_9(X)\},$
\par $\sigma_{SF_0}(T)=\{\lambda\in{\mathbb{C}}:  T-\lambda I\not\in SF_0(X)\},$
\par $\sigma_{lr}(T)=\{\lambda\in{\mathbb{C}}:  T-\lambda I\not\in M(X)\},$
\par The semi-regular spectrum: $\sigma_{se}(T)=\{\lambda\in{\mathbb{C}}:  T-\lambda I\not\in
D(X)\}$.

\vspace{3mm}

In [3], the semi-regular spectrum $\sigma_{se}(T)$ of $T$ is also
called as the regular spectrum and denoted by $ \sigma_{g}(T)$.

\vspace{3mm}

The spectra $\sigma_{BF}, \sigma_{BW}, \sigma_{rD}, \sigma_{lD},
\sigma_{SF_0}, \sigma_{lr}$ are closed and $\sigma_{SF_0},
\sigma_{lr}, \sigma_{se}$ are nonempty ([1], [3-5], [26-27]).

\vspace{3mm}

Let $D(\lambda, r)$ be the open disc centered at $\lambda\in
\mathbb{C}$ with radius $r>0$, $\overline{D}(\lambda, r)$ be the
corresponding closed disc. $T\in B(X)$ is said to have the Single
Valued Extension Property (SVEP, for short) ([28]) at $\lambda$ if
there exists $r
>0$ such that for each open subset $V\subseteq D(\lambda,r)$, the
constant function $f\equiv0$ is the only analytic solution of the
equation $$(T-\mu)f(\mu)=0, \forall \mu \in V.$$

For $T\in B(X)$, we denote $$S(T)=\{\lambda\in{\mathbb{C}}: T
\makebox{  fails to have the SVEP at }\lambda \}.$$ If
$S(T)=\emptyset$, then $T$ is said to have SVEP.

Henceforth, for $A\in B(X)$, $B\in B(Y)$ and $C\in B(Y, X)$, we put
$M_{C}=\left(
       \begin{array}{cc}
        A&C\\
         0&B\\

       \end{array}
     \right)$. It is clear that $M_C\in B(X\oplus Y)$.

\vspace{3mm}

For a given bounded linear operator $T\in B(X)$, as showed above, we
have defined 32 kinds spectra. Now, we are interesting in deciding
the spectra structure of operator matrix $M_{C}$. [28] and [16] have
told us that for spectra $\sigma, \sigma_e, \sigma_w, \sigma_a,
\sigma_{su}$ and $\sigma_{aw}$, we have

$$\sigma (M_C)\cup
(S(A^*)\cap S(B)) = \sigma (A)\cup \sigma (B),
\qquad\qquad\qquad\qquad\eqno(1)$$
$$\sigma_{e} (M_C)\cup (S(A^*)\cap
S(B)) = \sigma_{e} (A)\cup \sigma_{e}(B)\cup (S(A^*)\cap S(B)),
\eqno(2)$$
$$\sigma_{w}(M_C)\cup [S(A)\cap S(B^*)]\cup [S(A^*)\cap S(B)]=\qquad\qquad\qquad$$$$\sigma_{w}(A)\cup
\sigma_{w}(B)\cup [S(A)\cap S(B^*)]\cup [S(A^*)\cap S(B)],
\qquad\qquad\eqno(3)$$
$$\sigma_{su} (M_C)\cup  S(B)=\sigma_{su} (A)\cup
\sigma_{su}(B)\cup  S(B), \qquad\qquad\qquad\qquad\eqno(4)$$
 $$\sigma_{a} (M_C)\cup S(A^*)=\sigma_{a} (A)\cup \sigma_{a}(B)\cup
S(A^*), \qquad\qquad\qquad\qquad\eqno(5)$$
$$\sigma_{aw} (M_C)\cup (S(A)\cap S(B^*))\cup S(A^*) = \sigma_{aw}
(A)\cup \sigma_{aw}(B)\cup (S(A)\cap S(B^*))\cup S(A^*).\eqno(6)$$

\vspace{3mm}

In this paper, first, in Section 3, we decide another 20 kind
spectra structure of $M_C$, that is:

\vspace{3mm}

Spectra $\sigma_{D}$ and $ \sigma_{b}$ have equation (1) form,
spectra
    $\sigma_{des}, \sigma_{r}, \sigma_{rb},$ $ \sigma_{sb},\sigma_{re}$ and $\sigma_{SF-}$ have equation (4) form,
    spectra $\sigma_{l}, \sigma_{lb}, \sigma_{ab},$
$\sigma_{le}$ and $\sigma_{SF+}$ have equation (5) form,
    spectrum $\sigma_{lw}$ has equation (6) form. Moreover, if spectrum $\sigma_{*}=\sigma_{sw}$ or $\sigma_{rw}$, then it has the following form $$\sigma_{*}
(M_C)\cup (S(A)\cap S(B^*))\cup S(B) = \sigma_{*} (A)\cup
\sigma_{*}(B)\cup (S(A)\cap S(B^*))\cup S(B). \eqno(7)$$ If spectrum
$\sigma_{*}=\sigma_{K}, \sigma_{K_3}, \sigma_{SF0}$ or
$\sigma_{lr}$, then it has the following forms
$$\sigma_{*}(M_{C})\cup S(A^*)\cup S(A)\cup S(B)=\sigma_{*}(A)\cup \sigma_{*}(B)\cup S(A^*)\cup S(A) \cup
S(B)\eqno(8)$$ and $$\sigma_{*}(M_{C})\cup S(A)\cup S(B)\cup
S(B^*)=\sigma_{*}(A)\cup \sigma_{*}(B)\cup S(A) \cup S(B)\cup
S(B^*).\eqno(9)$$

\vspace{4mm}

On the other hand, Han and Lee in [19] studied the so-call
filling-in-hole problem of operator matrix, their result is:
$$\sigma_{}(A)\cup \sigma_{}(B)=\sigma_{}(M_C) \cup
W_{\sigma_{}}(A, B, C),\eqno(10) $$ where $W_{\sigma_{}}(A, B, C)$
is the union of some holes in $\sigma_{}(M_C)$ and $W_{\sigma_{}}(A,
B, C)\subseteq \sigma(A)\cap \sigma(B)$.

\vspace{3mm}

That is, the passage from $\sigma_{}(A)\cup \sigma_{}(B)$ to
$\sigma_{}(M_C)$ is the punching of some open sets in
$\sigma_{}(A)\cap \sigma_{}(B)$.

\vspace{3mm}

Moreover, in [12, 22, 31-32], the authors showed that for the
spectra $\sigma_{b}, \sigma_{w},
  \sigma_{e}$ and $ \sigma_{D}$, the equation (10) is also
  true.

\vspace{3mm}

In [12, 20], the authors showed that if spectrum
$\sigma_{*}=\sigma_{a}, \sigma_{SF+},\sigma_{SF-},\sigma_{ab}$ or
$\sigma_{sb}$, then $$\sigma_{*}(A)\cup \sigma_{*}(B)=(M_C) \cup
W_{\sigma_*}(A,B,C),\eqno(11)$$ where $W_{\sigma_*}(A, B, C)$ is
contained in the union of some holes in $\sigma_{*}(M_C)$.

\vspace{3mm}

Let $M$ be a compact subset of ${\mathbb{C}}$. The set
$$\eta{(M)}=\{w:\ |P(w)|\leq max\{|P(z)|:\ z\in M\}\mbox{ for every
polynomial}\, P\}$$ is called to be the polynomially convex hull of
$M$.

\vspace{3mm}

In [3], the authors proved that $$\eta(\sigma_{se}(A)\cup
\sigma_{se}(B))=\eta(\sigma_{se}(M_C)).\eqno(12)$$

Those spectra of $M_C$ which satisfy equation (10), (11), (12),
respectively, are called to have the filling-in-hole property,
generalized filling-in-hole property, convex filling-in-hole
property, respectively.

\vspace{3mm}

In Section 4, we continue to study the filling-in-hole problem of
another 18 kind spectra of $M_C$, we show that each spectrum of the
18 kind spectra of $M_C$ has one of the above 3 kind filling-in-hole
properties.

\vspace{3mm}

In Section 5, we present some interesting examples to show that some
conclusions about the spectra structure or filling-in-hole problem
of $M_C$ are not true.

%
%
%
%
\vskip 0.2in
 \large\section{Main Lemmas and Proofs}

\vskip 0.2in \noindent {\bf Lemma 2.1} ([31]). If $T \in B(X)$, then
the following statements are equivalent:
 \par (i). $T$ is Drazin invertible.
 \par (ii). $T=T_{1}\oplus T_{2},$  where  $T_{1}$ is invertible and  $T_{2}$ is nilpotent.
 \par (iii). There exists a nonnegative integer $k$ such that $des(T)=asc(T)=k<\infty.$
 \par (iv). $T^*$ is Drazin invertible.

\vskip 0.2in \noindent {\bf Lemma 2.2} ([31]). For $(A, B)\in
B(X)\times B(Y)$, if $M_{C}$ is Drazin invertible for some $C\in
B(Y,\,X)$, then
 \par (i). $des(B)<\infty$ and $asc(A)<\infty,$
 \par (ii). $des(A^*)<\infty$ and $asc(B^*)<\infty$.

\vskip 0.2in \noindent {\bf Lemma 2.3} ([1] Theorem 3.81). Let $T\in
B(X)$ and $des(T-\lambda_0)<\infty$. Then the following statements
are equivalent:
   \par (i). $T$ has the SVEP at    $\lambda_0$.
   \par(ii).  $asc(T-\lambda_0 )<\infty$.
   \par (iii). $\lambda_0 $ is a pole of the resolvent.
   \par(iv). $\lambda_0 $ is an isolated point of $\sigma(T)$.

\medskip\noindent {\bf Lemma 2.4} ([ 22, 31-32]). Let $(A, B)\in B(X)\times
B(Y)$ and $C\in B(Y,X)$. We have:

\par (i). if any two of operators $A, B$ and $M_C$ are invertible, then so is the third,
\par (ii). if any two of operators $A, B$ and $M_C$ are Fredholms, then so is the third,
\par (iii). if any two of operators $A, B$ and $M_C$ are Weyl, then so is the third,
\par (iv). if any two of operators $A, B$ and $M_C$ are Drazin invertible, then  so is the  third,
\par (v). if any two of operators $A, B$ and $M_C$ are Browder, then so is the third.

\vskip 0.2in \noindent {\bf Lemma 2.5} ([16]). Let $ (A, B)\in
B(X)\times B(Y )$ and $C\in B(Y, X)$. If $A$ has infinite ascent,
then $M_C$ has infinite ascent, if $B$ has infinite descent, then
$M_C$ has infinite descent.

 \vskip 0.2in

 \noindent {\bf Lemma 2.6}. Let $(A, B)\in B(X)\times B(Y)$ and $C\in B(Y, X)$. We have:
\par(i).  if $B$ is invertible, then $des(M_C)<\infty$ if and only if $des(A)<\infty$,
\par(ii).  if $ B=0$, then $des(M_C)<\infty$ if and only if $des(A)<\infty$,
\par(iii).  if $A$ is invertible, then $asc(M_C)<\infty$ if and only if $asc(B)<\infty$,
\par(iv).  if $A=0$, then $asc(M_C)<\infty$ if and only if $asc(B)<\infty$.

\medskip\par\noindent{\bf Proof.} (i). If $B$ is invertible, then $des(B)<\infty$, thus, it follows from $des(A)<\infty$ that
$des(M_C)<\infty$.  In fact, if $des(A)=p$ and $ des(B)=q$, then we
can claim that $R(M_C^{2n+1})=R(M_C^{2n})$ for each $C\in B(Y,X)$,
where  $n=max\{p,q\}$. For this, note that $R(M_C^{2n+1})\subseteq
R(M_C^{2n})$ is obvious, it is sufficient to show that
$R(M_C^{2n})\subseteq R(M_C^{2n+1})$.
\par  Suppose that  $u_0=(u_1,u_2)\in
R(M_C^{2n}).$  Then there exists $(x_0,y_0)\in X\oplus Y$ such that
$$(u_1,u_2) =M_C^{2n}(x_0,y_0)=(A^{2n}x_0+\sum
_{i=0}^{n-1}A^{2n-i-1}CB^{i}y_0+ \sum
_{i=n}^{2n-1}A^{2n-i-1}CB^{i}y_0,\ B^{2n}y_0),$$ where $A^0=I$ and
$B^0=I$. That $B^{2n}y_0=u_2$ is clear. In view of
$R(B^n)=R(B^{n+1}),$ there exists $y_1\in Y,$ such that
$$B^ny_0=B^{n+1}y_1.$$ Therefore,
\begin{eqnarray*}
              &u_1 &=A^{2n}x_0+\sum_{i=0}^{n-1}A^{2n-i-1}CB^{i}y_0+ \sum_{i=n}^{2n-1}A^{2n-i-1}CB^{i}y_0\\
   &      &=A^{2n}x_0+\sum_{i=0}^{n-1}A^{2n-i-1}CB^{i}y_0+ \sum_{i=n}^{2n-1}A^{2n-i-1}CB^{i+1}y_1\\
  &      &=A^{2n}x_0+\sum_{i=0}^{n-1}A^{2n-i-1}CB^{i}y_0+ \sum_{i=n+1}^{2n}A^{2n-i}CB^{i}y_1\\
   &      &=A^n(A^{n}x_0+\sum_{i=0}^{n-1}A^{n-i-1}CB^{i}y_0)+  \sum_{i=n+1}^{2n}A^{2n-i}CB^{i}y_1\\
   &      &=A^n(A^{n}x_0+\sum_{i=0}^{n-1}A^{n-i-1}CB^{i}y_0- \sum_{i=0}^{n}A^{n-i}CB^{i}y_1)+ \sum_{i=0}^{2n}A^{2n-i}CB^{i}y_1.
     \end{eqnarray*}
 Moreover, it follows from $des(A)=p\leq n<\infty$ that $R(A^n)=R(A^{n+1})=\cdots=R(A^{2n+1})$, hence, there exists $x_1\in
X$ such that $$A^n(A^{n}x_0+\sum_{i=0}^{n-1}A^{n-i-1}CB^{i}y_0-
\sum_{i=1}^{n-1}A^{n-i-1}CB^{i+1}y_1)=A^{2n+1}x_1.$$ Note that
$(u_1,u_2)=M_C^{2n+1}(x_1,y_1),$ so $R(M_C^{2n})\subseteq
R(M_C^{2n+1}).$

If $B$ is invertible and $des(M_C)<\infty$, now we consider two
cases to show that $des(A)<\infty$.

Case I.  If $R(M_C)=R(M_C^2)$, we claim that $R(A)=R(A^2)$. In fact,
$R(A)\supseteq R(A^2)$ is obvious. If $y\in R(A)$, there exists
$x\in X$ such that $y=Ax.$ Thus, $\left( \begin{array}{c} y\\0
\end{array} \right)=\left(
\begin{array}{cc}   A&C\\ 0&B\\ \end{array} \right) \left(
\begin{array}{c}   x\\ 0 \end{array} \right).$
That is $ \left( \begin{array}{c}   y\\0 \end{array} \right)\in
R(M_C)= R(M_C^2).$  Hence, there exists $ \left(
\begin{array}{c} x_1\\y_1\end{array} \right)\in X \oplus Y$ such
that $$ \left( \begin{array}{c}   y\\0 \end{array} \right) = {\left(
\begin{array}{cc}   A&C\\ 0&B\\ \end{array} \right) }^2 \left(
\begin{array}{c}    x_1\\y_1 \end{array} \right) = \left(
\begin{array}{cc}   A^2&AC+CB\\ 0&B^2\\ \end{array} \right) \left(
\begin{array}{c}    x_1\\y_1 \end{array} \right)=\left(
\begin{array}{c}   A^2 x_1+(AC+CB)y_1\\B^2y_1 \end{array} \right).$$
Moreover, since $B$ is invertible, it is easy to show that
$$y_1=0,y=A^2x_1.$$ Thus $y\in R(A^2),$ and so $R(A)\subseteq
R(A^2).$

 \par Case II.  If $1<des(M_C)=p<\infty$, put $M_{p}=M_C^p$.  Then
 $M_{p}= \left(\begin{array}{cc}   A^p&\sum _{i=1}^{p}A^{p-i}CB^{i-1}\\ 0&B^p\\ \end{array}
\right)$ and $R(M_{p})=R(M_{p}^2)$, where $A^0=I,B^0=I.$  By using
the same methods as in case I, we can prove that $R(A^p)\subseteq
R(A^{2p})$. So $R(A^p)=R(A^{2p})$. It follows from the conclusion
that $R(A^p)=R(A^{p+1}).$

\par Combining Case I with Case II, we complete the proof of (i).

\vskip 0.2in

\par (ii). It follows from the above argument methods that we only need to show
that if $B=0$ and $R(M_C)=R(M_C^2)$, then $R(A)=R(A^2)$. In fact, if
$y\in R(A)$, then there exists $x\in X$ such that $y=Ax,$ so $
\left( \begin{array}{c}   y\\0 \end{array} \right) = \left(
\begin{array}{cc}   A&C\\ 0&0\\ \end{array} \right) \left(
\begin{array}{c}   x\\ 0 \end{array} \right)$. Since
 $ \left( \begin{array}{c}   y\\0 \end{array} \right)\in
{R(M_C)= R(M_C^3)},$  there exists $\left(
\begin{array}{c} x_1\\y_1\end{array} \right)\in X \oplus Y$ such
that $$ \left( \begin{array}{c}   y\\0 \end{array} \right) = {\left(
\begin{array}{cc}   A&C\\ 0&0\\ \end{array} \right) }^3 \left(
\begin{array}{c}    x_1\\y_1 \end{array} \right)  =\left(
\begin{array}{cc}   A^2&0\\ 0&0\\ \end{array} \right) \left(
\begin{array}{c}   A x_1+Cy_1\\0 \end{array} \right).$$
Thus we have $y=A^2(A x_1+Cy_1)$, so $y\in R(A^2), $ this showed
that $R(A)\subseteq R(A^2)$. Note that $R(A^2)\subseteq R(A)$ is
obvious, therefore $R(A)= R(A^2).$

\vskip 0.2in

(iii). If $A$ is invertible, then $asc(A)<\infty$, thus, it follows
easily from $asc(B)<\infty$ that $asc(M_C)<\infty$ (see Lemma 2.2 in
[11]).

If $A$ is invertible and $asc(M_C)<\infty$, now we consider two
cases to show that $asc(B)<\infty$.

\par Case I. $N(M_C)=N(M_C^2)$. We claim that $N(B)=N(B^2).$
 If $y\in N(B^2)$, then
 $B^2y=0$. Note that $A$ is invertible, it can be proved that
 $$M_C^2\left( \begin{array}{c} -A^{-2} (AC+CB)  y\\y \end{array} \right)=\left( \begin{array}{c} 0\\B^2y \end{array} \right)
 =\left( \begin{array}{c} 0\\0 \end{array} \right).$$
 Thus, $$\left( \begin{array}{c} -A^{-2} (AC+CB)  y\\y \end{array} \right)\in  N(M_C^2)=N(M_C).  $$
Therefore  $$ \left( \begin{array}{cc}   A&C\\ 0&B\\ \end{array}
\right)\left( \begin{array}{c} -A^{-2} (AC+CB)  y\\y \end{array}
\right)= \left( \begin{array}{c}  0\\0 \end{array} \right),
$$ which implies $By=0$, thus $N(B^2)\subseteq N(B)$. Note that $N(B)\subseteq N(B^2) $ is
obvious, so $N(B^2)=N(B)$.

 \par Case II.
  If $1<asc(M_C)=p<\infty$, put $M_{p}=M_C^p$.  Then
 $$M_{p}= \left(
\begin{array}{cc}   A^p&\sum _{i=1}^{p}A^{p-i}CB^{i-1}\\ 0&B^p\\ \end{array}
\right)\,\,\makebox{ and}\,\, N(M_{p})=N(M_{p}^2),$$ where
$A^0=I,B^0=I.$ By the same methods as in case I, we can prove that
$N(B^p)\subseteq N(B^{2p})$ and so it is easy to obtain that
$N(B^p)= N(B^{p+1}).$
\par Combining case I with case II,  we prove (iii).

\vspace{3mm}

\par (iv). We only prove that
 if $A=0$ and $N(M_C)=N(M_C^2)$, then $N(B)=N(B^2).$
Let $y\in N(B^2)$. Then
 $B^2y=0$.  So we have
 $$M_C^3\left( \begin{array}{c} 0\\y \end{array} \right)
 =\left( \begin{array}{cc}   0&CB^2\\ 0&B^3\\ \end{array} \right)\left( \begin{array}{c} 0\\y \end{array} \right)
 =\left( \begin{array}{c} 0\\0 \end{array} \right),$$
 which implies that $$\left( \begin{array}{c}0\\y \end{array} \right)\in  N(M_C^3)=N(M_C).  $$
Thus $$ \left( \begin{array}{cc}   0&C\\ 0&B\\ \end{array}
\right)\left( \begin{array}{c} 0\\y \end{array} \right)= \left(
\begin{array}{c}  Cy\\By \end{array} \right)= \left(
\begin{array}{c}  0\\0 \end{array} \right).$$ This showes that
$By=0,$ so $N(B^2)\subseteq N(B)$. It follows easily from the
conclusion that $N(B)= N(B^2).$ The lemma is proved.

\vspace{3mm}
 The following lemma is important and it is widely used in the
 proofs of our main theorems in Section 3 and Section 4.

 \vspace{3mm}

\noindent {\bf Lemma 2.7.} For $(A, B)\in B(X)\times B(Y)$ and $C\in
B(Y,X)$, we have:
 \par (i).  if $A$ is Drazin invertible, then $asc(M_C)<\infty$ ($des(M_C)<\infty$) if and only if  $asc(B)<\infty$ ($des(B)<\infty$),
 \par (ii). if $B$ is Drazin invertible, then $des(M_C)<\infty$ ($asc(M_C)<\infty$ ) if and only if  $des(A)<\infty$ ($asc(A)<\infty$),
 \par (iii). if $A$ is Browder operator, then $M_C$ is left (right, upper, lower) semi-Browder operator if and only if  $B$ is left (right, upper, lower) semi-Browder operator,
\par  (iv).  if $B$ is Browder operator, then  $M_C$ is left (right, upper, lower) semi-Browder operator if and only if $A$ is left  (right, upper, lower)  semi-Browder operator,
 \par (v). if $A$ is Fredholm operator, then $M_C$ is left  (right, upper, lower)  semi-Fredholm operator if and only if  $B$ is left  (right, upper, lower)  semi-Fredholm operator,
\par  (vi). if $B$ is Fredholm operator, then  $M_C$ is left  (right, upper, lower)  semi-Fredholm operator if and only if $A$ is left  (right, upper, lower)  semi-Fredholm operator,
\par  (vii). if $A$ is Weyl operator, then $M_C$ is left  (right, upper, lower)  semi-Weyl operator if and only if  $B$ is left  (right, upper, lower)  Weyl operator,
\par  (viii). if $B$ is Weyl operator, then  $M_C$ is left  (right, upper, lower)  semi-Weyl operator if and only if $A$ is left  (right, upper, lower)  semi-Weyl operator,
\par  (viiii). if $A$ is invertible, then $M_C$ is left (right) invertible if and only if  $B$ is left (right)invertible,
\par (vv). if $B$ is invertible, then $M_C$ is left (right) invertible if and only if $A$ is left (right) invertible,
\par (vvi). if $A$ is invertible, then $M_C$ bounded below if and only if  $B$  bounded below,
\par (vvii). if $B$ is invertible, then $M_C$ is  surjective if and only if $A$ is  surjective.

\medskip\par\noindent{\bf Proof.} That from (v) to (vv) are well known and from (i) to (iv) are new. Here we only prove (i).

\par  If $A$ is Drazin invertible and $i(A)=k<\infty$, then $$A=\left( \begin{array}{cc}  A_1&0\\ 0&A_2\\\end{array}\right):
  R(A^k)\oplus  N(A^k)\longrightarrow R(A^k)\oplus  N(A^k),$$
where $A_1$ is invertible and $A_2^k=0.$
 Thus,
  $$M_C=\left( \begin{array}{ccc}   A_1&0&C_1\\ 0&A_2&C_2\\ 0&0&B\\ \end{array}\right):
 R(A^k)\oplus  N(A^k) \oplus Y\longrightarrow  R(A^k)\oplus  N(A^k)\oplus Y.$$
So
  $$M_C^k=\left( \begin{array}{ccc}   A_1^k&0&\sum _{i=1}^{k}A_1^{k-i}C_1B^{i-1}\\ 0&0&\sum _{i=1}^{k}A_2^{k-i}C_2B^{i-1}\\ 0&0&B^k\\ \end{array}\right):
  R(A^k)\oplus  N(A^k) \oplus Y\longrightarrow   R(A^k)\oplus  N(A^k)\oplus Y.$$

It follows from Lemma 2.6 that
\begin{eqnarray*}
& asc(M_C)<\infty &\Longleftrightarrow asc(M_C^k)<\infty\\
& &\Longleftrightarrow asc(\left( \begin{array}{cc}   0&\sum _{i=1}^{k}A_2^{k-i}C_2B^{i-1}\\ 0&B^k\\ \end{array}\right))<\infty \\
& &\Longleftrightarrow asc(B^k)<\infty\\
& &\Longleftrightarrow asc(B)<\infty.
\end{eqnarray*}
 Moreover, if $A$ is Drazin
invertible, then it follows from Lemma 2.5 and the proof of Lemma
2.6 that
$$des(M_C)<\infty\Leftrightarrow des(B)<\infty.$$ Thus (i) is proved.

\vskip 0.2in \noindent {\bf Lemma 2.8} ([1] Theorem 3.4). For $T\in
B(X)$, the following  properties hold:
\par (i). if $asc(T)<\infty$, then $ \alpha(T)\leq\beta(T)$,
\par(ii). if $des(T)<\infty$, then $ \beta(T)\leq\alpha(T)$,
\par(iii). if $asc(T)=des(T)<\infty$, then $\beta(T)=\alpha(T)$,
\par(iv). if $\beta(T)=\alpha(T)<\infty$ and either $asc(T)$ or $des(T)$ is finite, then $ asc(T)=des(T)$.

\vskip 0.2in \noindent {\bf Lemma 2.9} ([1] Corollary 3.19).
 Let $T\in B(X)$ and  $\lambda_0-T\in \Phi_\pm(X)$, where $\Phi_\pm(X)=\Phi_+(X)\cup\Phi_-(X)$.
  We have:

  \par (i). if $T$ has the SVEP at $\lambda_0$, then ind $(\lambda_0I-T)\leq0$,
  \par(ii). if $T^*$ has the SVEP at $\lambda_0$, then ind $(\lambda_0I-T)\geq0$.

\vskip 0.2in \noindent {\bf Lemma 2.10} ([31]). If each neighborhood
of $\lambda $ contains a point that is not an eigenvalue of operator
$T$ and $\lambda-T $ has a  finite descent, then $\lambda-T $ is
Drazin invertible.

%
%
%
%
\large\section{Spectra structure of operator matrix $M_C$ }

The following theorems in this section tell us 20 kind spectra
structure of $M_C$, that is:

\vskip 0.2in \noindent {\bf Theorem 3.1.} $$\sigma _D(M_C)\cup
(S(A^*)\cap S(B)) = \sigma _D(A)\cup \sigma _D(B).$$

\par\noindent{\bf Proof.} First, Theorem 2.9 in [31] told us that $\sigma _D(M_C)\subseteq \sigma _D(A)\cup
\sigma_D(B)$.
 Note that  $$S(A^*)\cap
S(B)\subseteq int( \sigma (A^*))\cap int (\sigma (B))=int( \sigma
(A))\cap int (\sigma (B))\subseteq \sigma _D(A)\cup \sigma _D(B).$$
It follows  that $$\sigma _D(M_C)\cup (S(A^*)\cap S(B))\subseteq
\sigma _D(A)\cup \sigma _D(B).$$

\par On the other hand, if $\lambda \in (\sigma _D(A)\cup \sigma
     _D(B))\setminus\sigma _D(M_C)$, then $M_C-\lambda$ is Drazin
     invertible. It follows from Lemma 2.2 that $$des(B-\lambda)<\infty \,\,\makebox{
     and}
     \,\,des(A^*-\lambda)<\infty.$$

Now we claim that $\lambda\in S(A ^* ) \cap S(B)$. In fact, if
$\lambda\not \in S (A^*) \cap S (B)$. There are two cases:

If $\lambda\not \in S (A^*)$, note that $des(A^*-\lambda)<\infty$,
it follows from Lemmas 2.1 and Lemma 2.3  that
$asc(A^*-\lambda)<\infty$ and $ A -\lambda$ is Drazin invertible.
Furthermore, Lemma 2.4 tells us that $B-\lambda$ is also Drazin
invertible, this is a contradiction.

Similarly, we can prove that $\lambda\not \in S(B)$ is also
impossible. Thus, we have $$\sigma _D(A)\cup \sigma _D(B) \subseteq
\sigma _D(M_C)\cup (S(A^*)\cap S(B)).$$ The theorem is proved.

\vskip 0.2in \noindent {\bf Theorem 3.2.} $$ \sigma _b(M_C)\cup
(S(A^*)\cap S(B)) = \sigma _b(A)\cup \sigma _b(B).$$

\par\noindent{\bf Proof.}  It follows from Lemma 2.4 that $ \sigma _b(M_C)\subseteq\sigma _b(A)\cup \sigma
 _b(B).$  Moreover, it is easy to know that $\sigma_{D}(A)\cup \sigma _{D}(B)\subseteq \sigma_{b}(A)\cup \sigma _{b}(B)$,
 thus, it follows from the proof of Theorem 3.1 that $S(A^*)\cap S(B) \subseteq \sigma _b(A)\cup \sigma _b(B).$
  Hence $$ \sigma _b(M_C)\cup (S(A^*)\cap S(B))\subseteq \sigma _b(A)\cup \sigma_b(B).$$

 For the  contrary inclusion, it is sufficient to prove that
  $$(\sigma _b(A)\cup \sigma_b(B))\setminus\sigma _b(M_C) \subseteq (S(A^*)\cap S(B)).$$

\par Let $\lambda \in (\sigma _b(A)\cup \sigma
     _b(B))\setminus\sigma _b(M_C)$. Then $M_C-\lambda$ is a Drazin
     invertible Fredholm operator, and so $A-\lambda\in \Phi_+(X), B-\lambda\in \Phi_-(Y) $.
     Moreover, It follows from Lemma 2.2 that $$des(B-\lambda)<\infty \,\,\makebox{ and}
     \,\,asc(A-\lambda)<\infty.$$

Now, we show that $\lambda\in S(A ^* ) \cap S(B)$. In fact, if
$\lambda\not \in S (A^*) \cap S (B)$, then when $\lambda\not \in S
(A^*)$, we can prove as in Theorem 3.1 that $A-\lambda$ is Drazin
invertible. Note that $A-\lambda\in \Phi_+(X)$, it is easy to show
from Lemma 2.8 that $A-\lambda\in \Phi_b(X).$ Thus, by Lemma 2.4 we
get  that $B-\lambda\in \Phi_b(Y)$, and hence $\lambda \not\in
\sigma _b(A)\cup \sigma
     _b(B)$, which is a contradiction. So, $\lambda\in S (A^*)$.

Similarly, we can show that $\lambda\in S (B)$. The theorem is
proved.

\vskip 0.2in \noindent {\bf Theorem 3.3.} If $
\sigma_{*}=\sigma_{des}, \sigma_{r},
    \sigma_{rb}, \sigma_{sb}, \sigma_{re}$ or $\sigma_{SF-}$, then
 $$\sigma_{*}(M_{C})\cup S{(B)}=\sigma_{*}(A)\cup
\sigma_{*}(B) \cup S{(B)}.$$

  \par\noindent{\bf Proof.} Observe that
$$ M_{C} = \left( \begin{array}{cc} I&0\\ 0&B\\ \end{array} \right)
     \left(\begin{array}{cc} I&C\\ 0&I\\ \end{array} \right)
     \left(   \begin{array}{cc}  A&0\\  0&I\\  \end{array} \right).$$

 \noindent It is easy to prove that $$\sigma_{*}(B)\subseteq \sigma_{*}(M_{C})\subseteq
\sigma_{*}(A)\cup \sigma_{*}(B)$$ and
$$\sigma_{*}(M_{C})\cup S{(B)}\subseteq \sigma_{*}(A)\cup
\sigma_{*}(B)\cup S{(B)}.$$
 \noindent for each $\sigma_{*}=\sigma_{des}, \sigma_{r}, \sigma_{rb}, \sigma_{sb}, \sigma_{re}$
 or $\sigma_{SF-}$.  So, in order to prove the theorem,
 we only need to prove that
  $$(\sigma _{*}(A)\cup \sigma
_{*}(B))\setminus\sigma _{*}(M_C) \subseteq S(B).$$
 If $\lambda \in (\sigma _{*}(A)\cup \sigma
     _{*}(B))\setminus\sigma _{*}(M_C)$ and $\sigma_{*}=\sigma_{des}, \sigma_{r},
    \sigma_{rb}, \sigma_{sb}, \sigma_{re}$ or $\sigma_{SF-}$, we show that $\lambda\in S(B
    ).$ In fact, if $\lambda\not \in S (B)$, we consider the following four
cases:

Case (I). Let $\sigma_{*}=\sigma_{des}$. Since $\lambda \in (\sigma
_{des}(A)\cup \sigma
     _{des}(B))\setminus\sigma _{des}(M_C)$ and $\lambda\not \in S
     (B)$, it follows from Lemmas 2.5 and 2.3 that $des(B- \lambda )<\infty$ and $ asc(B-
     \lambda)<\infty$, so $B-
     \lambda$ is Drazin invertible.  Thus, Lemma 2.7 tells us that $des(A- \lambda)<\infty,$ this is a contradiction.
     Therefore,  it follows that $(\sigma _{des}(A)\cup \sigma_{des}(B))\setminus\sigma _{des}(M_C) \subseteq S(B)$.
\par
 Case (II). Let $\sigma_{*}=\sigma_{r}$. Note that
$\lambda \in (\sigma _{r}(A)\cup \sigma
     _{r}(B))\setminus\sigma _{r}(M_C)$, so $B- \lambda \in G
     _{r}(Y).$ That is $B- \lambda \in\Phi
     _{rb}(Y)$  and $R(B- \lambda)=Y$. Since $\lambda\not \in S
     (B)$, it follows from Lemma 2.9 that ind$(B- \lambda)=\alpha(B- \lambda)-\beta(B-
     \lambda)\leq0$. Thus $\alpha(B- \lambda)=0$,
     and so $B- \lambda \in \Phi_0(Y)$.  Also since $R(B- \lambda)=Y$, $B- \lambda$ is invertible.
     It follows from Lemma 2.7 that $A- \lambda \in G_r(X)$, this contradicts $\lambda \in (\sigma _{r}(A)\cup \sigma
     _{r}(B))$. Therefore $(\sigma _{r}(A)\cup \sigma
_{r}(B))\setminus\sigma _{r}(M_C) \subseteq S(B)$.

Case (III). Let $\sigma_{*}=\sigma_{rb}$. Since $\lambda \in (\sigma
_{rb}(A)\cup \sigma
     _{rb}(B))\setminus\sigma _{rb}(M_C)$,   $B- \lambda \in\Phi
     _{rb}(Y).$ In particular,  $B- \lambda \in\Phi
     _{r}(Y)$, ind$(B- \lambda)\geq0$. Note that $\lambda\not \in S
     (B)$, it follows from Lemma 2.9 that ind$(B-
     \lambda)\leq 0$. Thus $B- \lambda \in \Phi
     _{0}(Y)$ and  $des(B- \lambda)<\infty$,  and then Lemma 2.8 shows  us that $B- \lambda \in \Phi
     _{b}(Y)$. Lemma 2.7 tells us that $A- \lambda \in \Phi
     _{rb}(X),$ this is a contradiction. So we get that
     $(\sigma _{rb}(A)\cup \sigma
_{rb}(B))\setminus\sigma _{rb}(M_C) \subseteq S(B)$.

Case (IV). Let $\sigma_{*}=\sigma_{re}$.
 Since $\lambda \in (\sigma _{re}(A)\cup \sigma
     _{re}(B))\setminus\sigma _{re}(M_C)$, it is easy to show that $B- \lambda \in\Phi
     _{r}(Y)$. Note that  $\lambda\not \in S
     (B)$, it follows from Lemma 2.9 that  ind$(B-
     \lambda)\leq 0$, so $\alpha(B- \lambda) \leq\beta(B-
     \lambda)<\infty.$
     Thus $B- \lambda \in \Phi
     _{}(X)$. By using Lemma 2.7, we have $A- \lambda \in \Phi
     _{r}(X),$ this is a contradiction and so we can prove that $(\sigma _{re}(A)\cup \sigma
_{re}(B))\setminus\sigma _{re}(M_C) \subseteq S(B)$. For
$\sigma_{*}=\sigma_{sb}$ or $\sigma_{SF-}$, the proof methods are
similar, we omit them.

\vspace{3mm}

\par Similarly, we can prove also the following theorem:

\vskip 0.2in \noindent {\bf Theorem 3.4.} If $
\sigma_{*}=\sigma_{l},
    \sigma_{lb}, \sigma_{ab}, \sigma_{le}$ or $\sigma_{SF+}$, then
     $$\sigma_{*}(M_{C})\cup S{(A^*)}=\sigma_{*}(A)\cup \sigma_{*}(B) \cup
S{(A^*)}.$$

\vskip 0.2in \noindent {\bf Theorem 3.5.}  $$\sigma_{lw} (M_C)\cup
(S(A)\cap S(B^*))\cup S(A^*) = \sigma_{lw} (A)\cup
\sigma_{lw}(B)\cup (S(A)\cap S(B^*))\cup S(A^*).$$

\par\noindent{\bf Proof.}  Note that $\sigma_{lw} (M_C)\subseteq
\sigma_{lw} (A)\cup \sigma_{lw}(B),$  it is obvious that
$$\sigma_{lw} (M_C)\cup (S(A)\cap S(B^*))\cup S(A^*)\subseteq
\sigma_{lw} (A)\cup \sigma_{lw}(B)\cup (S(A)\cap S(B^*))\cup
S(A^*).$$

 For the contrary inclusion, let $\lambda\not\in\sigma_{lw} (M_C)\cup (S(A)\cap S(B^*))\cup
S(A^*).$  Then $M_C-\lambda$ is left semi-Weyl operator, and it is
easy to prove that $A-\lambda\in\Phi_l(X)$. From the assumption we
also get that either $A$ and $A^*$ or $A^*$ and $B^*$ have the SVEP
at $\lambda$.
 If $A$ and $A^*$  have the SVEP at $\lambda$,  it follows from Lemma 2.9 that
 $A-\lambda\in \Phi_0(X)$. By using Lemma 2.7, it follows that  $B-\lambda\in
 \Phi_{l}( Y)$  and  ind$(M_C-\lambda)=$
 ind$(A-\lambda)+$ind$(B-\lambda)\leq0$.
 Thus  $B-\lambda\in \Phi_{l}(Y)$ with ind$(B-\lambda)=$ind$(M_C-\lambda)$-ind$(A-\lambda)$=ind$(M_C-\lambda)\leq0.$
 Hence $\lambda\not\in \sigma_{lw} (A)\cup \sigma_{lw}(B).$
 On the other hand, if $A^*$ and $B^*$ have the SVEP at $\lambda$,
it is  obvious that  $M_C^*$ have the SVEP at $\lambda$. It follows
from Lemma 2.9  that $M_C -\lambda\in \Phi_0(X\oplus Y)$.  Thus, it
is easy to show that  $A -\lambda\in \Phi_l(X)$ and $B -\lambda\in
\Phi_r(Y)$.  Also, note that $B^*$ and $A^*$ have the SVEP at
$\lambda$, it follows from Lemma 2.9 that ind$(A-\lambda)\geq0$ and
ind$(B-\lambda)\geq0$.  Hence it follows from $A -\lambda\in
\Phi_l(X)$ proved above that $A-\lambda\in
 \Phi_{}(X)$, so Lemma 2.4 tells us that
$B-\lambda\in \Phi_{}(Y)$. Moreover, in view of
$0=$ind$(M_C-\lambda)=$ind$(A-\lambda)+$ind$(B-\lambda)$,
ind$(A-\lambda)\geq0$ and ind$(B-\lambda)\geq0$, it is clear that
ind$(A-\lambda)$=ind$(B-\lambda)=0$. Thus $A-\lambda$ and
$B-\lambda$ are both Weyl operators, which implies that
$\lambda\not\in \sigma_{lw} (A)\cup \sigma_{lw}(B).$  It follows
that $\lambda\not\in \sigma_{lw} (A)\cup \sigma_{lw}(B)$ when
$\lambda\not\in\sigma_{lw} (M_C)\cup (S(A)\cap S(B^*))\cup S(A^*).$
Thus, the contrary inclusion is clear. The theorem is proved.

\vspace{3mm}

\par  Similarly, we can prove also the following theorem:

 \vskip 0.2in \noindent {\bf Theorem 3.6.} If $\sigma_{*}=\sigma_{sw}$ or $\sigma_{rw}$, then $$\sigma_{*}
(M_C)\cup (S(A)\cap S(B^*))\cup S(B) = \sigma_{*} (A)\cup
\sigma_{*}(B)\cup (S(A)\cap S(B^*))\cup S(B).$$

\vskip 0.2in \noindent {\bf Theorem 3.7.} If $\sigma_{*}=\sigma_{K},
\sigma_{K_3}, \sigma_{SF0}$ or
 $\sigma_{lr}$, then $$\sigma_{*}(M_{C})\cup S(A^*)\cup S(A)\cup S(B)=\sigma_{*}(A)\cup \sigma_{*}(B)\cup S(A^*)\cup
S(A) \cup S(B)$$ and
$$\sigma_{*}(M_{C})\cup S(A^*)\cup S(B)\cup S(B^*)=\sigma_{*}(A)\cup \sigma_{*}(B)\cup S(A^*) \cup S(B)\cup
S(B^*).$$

\par\noindent{\bf Proof.} We only prove when $\sigma_{*}=\sigma_{K}$, equation
(8) holds. First, we prove that
$$\sigma_{K}(M_{C})\cup S(A^*)\cup S(A)\cup S(B)\subseteq \sigma_{K}(A)\cup \sigma_{K}(B)\cup S(A^*)\cup S(A) \cup
S(B).$$ In fact, let $\lambda\not\in \sigma_{K}(A)\cup
\sigma_{K}(B)\cup S(A^*)\cup S(A) \cup S(B).$ Then  $A-\lambda$  is
a semi-Fredholm operator with $A^*$, $A$   have the SVEP at
$\lambda$, this implies $A-\lambda$ is a Weyl operator. Note that
$B-\lambda$  is also a semi-Fredholm operator, it follows from Lemma
2.7 that $M_C-\lambda$ is a semi-Fredholm operator. Thus
$\lambda\not\in\sigma_{K}(M_{C})\cup S(A^*)\cup S(A)\cup S(B).$ So
it is clear that $\sigma_{K}(M_{C})\cup S(A^*)\cup S(A)\cup
S(B)\subseteq \sigma_{K}(A)\cup \sigma_{K}(B)\cup S(A^*)\cup S(A)
\cup S(B).$

 For the contrary inclusion, let $\lambda\not\in\sigma_{K}(M_{C})\cup
S(A^*)\cup S(A)\cup S(B).$  Then $M_C-\lambda$ is a semi-Fredholm
operator, that is, $M_C-\lambda$ is either a upper semi-Fredholm
operator or  a lower semi-Fredholm operator. If
$M_C-\lambda\in\Phi_+(X\oplus Y)$, then it is easy to prove that
$A-\lambda \in\Phi_+(X)$.  Since $A^*$ has the SVEP at $\lambda$, it
follows from lemma 2.9  that $A-\lambda \in\Phi (X)$.  Hence it
follows Lemma 2.7 that $B-\lambda \in\Phi_+ (Y)$. Thus
$\lambda\not\in \sigma_{K}(A)\cup \sigma_{K}(B)$, and so
$\lambda\not\in \sigma_{K}(A)\cup \sigma_{K}(B)\cup S(A^*)\cup S(A)
\cup S(B).$ On the other hand, if $M_C-\lambda\in\Phi_-(X\oplus Y)$,
then $B-\lambda \in\Phi_-(Y)$. Since $B$ has the SVEP at $\lambda$,
it follows from lemma 2.9  that $B-\lambda \in\Phi (Y)$. Similar to
the above arguments, we can also obtain that $\lambda\not\in
\sigma_{K}(A)\cup \sigma_{K}(B)\cup S(A^*)\cup S(A) \cup S(B)$. This
proves equation (8).

\vspace{3mm}

We are interesting in the following question, it is perhaps
difficult:

\vspace{3mm}

\medskip\par {\bf Open question 3.8.} Do other spectra of $M_C$ have the equations (1) to (9) forms ?

%
%
%
%
\large\section{Filling-in-hole Problem of $M_C$ }
In Section 1, we pointed out that if spectrum $\sigma_{*}=\sigma,
\sigma_{b}, \sigma_{w}, \sigma_{e}$ or $\sigma_{D}$, then
$$\sigma_{*}(A)\cup \sigma_{*}(B)=\sigma_{*}(M_C) \cup
W_{\sigma_{*}}(A, B, C),$$ where $W_{\sigma_{*}}(A, B, C)$ is the
union of some holes in $\sigma_{*}(M_C)$ and $W_{\sigma_{*}}(A, B,
C)\subseteq \sigma_{*}(A)\cap \sigma_{*}(B)$. The following theorem
shows the relationship among $W_{\sigma_{}}(A,B,C),
W_{\sigma_{b}}(A,B,C)$, $ W_{\sigma_{D}}(A,B,C)$ and
$W_{\sigma_{w}}(A,B,C)$:

\vskip 0.2in\noindent

 \noindent {\bf Theorem 4.1.} For $(A, B)\in B(X)\times B(Y)$ and $C\in B(Y, X)$, we have

 \vspace{2mm}

 (i). $ W_{\sigma_{}}(A,B,C)\subseteq W_{\sigma_{b}}(A,B,C)\subseteq W_{\sigma_{D}}(A,B,C),$

 \vspace{2mm}

 (ii). $ W_{\sigma_{b}}(A,B,C)\subseteq W_{\sigma_{w}}(A,B,C).$

 \vspace{3mm}

In particular,  the following states are equivalent:

\vspace{3mm}

 \par (a). $W_{\sigma_{}}(A,B,C)=\emptyset,$

 \vspace{2mm}

 \par (b). $W_{\sigma_{b}}(A,B,C)=\emptyset,$

 \vspace{2mm}

  \par (c). $W_{\sigma_{D}}(A,B,C)=\emptyset.$

\par\vskip 0.2in \noindent{\bf Proof.} (i). Let
$\lambda \in W _{\sigma_{}}(A,B,C)$. It follow from equation (1)
that $\lambda \in (\sigma (A)\cup \sigma_{}(B))\setminus
\sigma_{}(M_{C})$. Thus $A-\lambda$ is left invertible and
$B-\lambda$ is right invertible. That $\lambda \not\in
\sigma_{b}(M_{C})$ is obvious. Now we claim that $\lambda \in
\sigma_{b}(A)\cup \sigma_{b}(B)$. If not, by Lemma 2.4, we have that
both $A-\lambda$ and $B-\lambda$ are Browder operators. This implies
that $\lambda \in \Phi_0 (A)\cap \Phi_{0}(B)$. Moreover, since
$A-\lambda$ is left invertible and $B-\lambda$ is right invertible,
$A-\lambda$ and $B-\lambda$ are invertible, this contradicts
$\lambda \in \sigma (A)\cup \sigma_{}(B)$, so $\lambda \in \sigma_b
(A)\cup \sigma_{b}(B)$, thus $ W_{\sigma_{}}(A,B,C)\subseteq
W_{\sigma_{b}}(A,B,C)$ is proved.

\vspace{2mm}

 For the inclusion $W_{\sigma_{b}}(A,B,C)\subseteq
 W_{\sigma_{D}}(A,B,C)$, note that $\sigma_{b}(M_{C}) \supseteq
 \sigma_{D}(M_{C})$, so it is sufficient to show that if $\lambda \in (\sigma_{b}(A)\cup\sigma_{b}(B))\setminus
\sigma_{b}(M_{C})$, then
   $\lambda \in (\sigma_{D}(A)\cup  \sigma_{D}(B))$. Let $\lambda \in (\sigma_{b}(A)\cup \sigma_{b}(B))\setminus
\sigma_{b}(M_{C})$. Then $M_{C}-\lambda\in \Phi_b (X\oplus Y)$, so
$A-\lambda \in \Phi_l( X)$ and $B-\lambda \in \Phi_r( Y)$. We claim
$\lambda \in \sigma_D (A)\cup  \sigma_{D}(B)$. If not, it follows
from Lemma 2.4 that $\lambda \in \rho_D (A)\cap \rho_{D}(B)$. Since
$B-\lambda \in \Phi_r( Y)$,
 Lemma 2.1 and  Lemma 2.8 tell us that $B-\lambda$ is a Fredholm operator.
  This implies that $B-\lambda$ is a Drazin invertible Fredholm operator, that is, $B-\lambda$ is a Browder operator.
 By Lemma 2.4, we get that $A-\lambda$ is also a Browder operator. Thus $\lambda \in (\Phi _{b }(A)\cap  \Phi_{b}(B))$, which contradicts
 with the assumption that $\lambda \in (\sigma _{b }(A)\cup  \sigma_{b}(B))$, so we have $\lambda \in \sigma_D (A)\cup  \sigma_{D}(B)$, thus $W_{\sigma_{b}}(A,B,C)\subseteq W_{\sigma_{D}}(A,B,C)$ is also proved.

\vspace{3mm}

 (ii). To prove  $W_{\sigma_{b}}(A,B,C)\subseteq W_{\sigma_{w}}(A,B,C)$, by Lemma 2.4 and the fact that $\sigma_{w}(M_{C})\subseteq\sigma_{b}(M_{C})$,
 it sufficient to show that if $\lambda \in (\sigma _{b }(A)\cup  \sigma_{b}(B))\setminus
\sigma_{b}(M_{C})$, then
   $\lambda \in \sigma_{w}(A)\cup  \sigma_{w}(B)$. Indeed, if $\lambda \in (\sigma_ {b} (A)\cup  \sigma_{b}(B))\setminus \sigma_{b}(M_{C})$,
   then $M_{C}-\lambda\in \Phi_b(X\oplus Y)$.
It follows Lemma 2.2 that  there exist some nonnegative integer $k$
and $l$ such that $asc(A-\lambda)=k<\infty$ and
$des(B-\lambda)=l<\infty$. Now we claim $\lambda \in \sigma_w
(A)\cup  \sigma_{w}(B)$. Otherwise, $\lambda \in \Phi_0 (A)\cap
\Phi_{0}(B)$. Moreover, Since $asc(A-\lambda)=k<\infty$, by Lemma
2.8 we have $des(A-\lambda)<\infty.$  That is  $A-\lambda$ is a
Drazin invertible Fredholm operator, so $A-\lambda\in \Phi_b(X).$
Using Lemma 2.4, we get that $B-\lambda\in \Phi _b(Y)$ and so
$\lambda\in \Phi _b(A)\cap \Phi _b(B)$, this contradicts with the
assumption that $\lambda \in (\sigma _{b }(A)\cup \sigma_{b}(B))$.
Thus we have $\lambda \in \sigma_w (A)\cup \sigma_{w}(B)$ and (ii)
is proved.

\vspace{3mm}

Finally, it follows from Corollary 2.12 in [31] that if $\
W_{\sigma_{}}(A,B,C)=\emptyset$, then $
W_{\sigma_{D}}(A,B,C)=\emptyset.$ This completed the proof of
theorem.

\vspace{3mm}

The following result which generalizes Lemma 3.2 of [12] is useful
in studying the filling-in-hole problem of $M_C$.

\vskip 0.2in\noindent {\bf Proposition 4.2.} For each $T\in B(X)$,
we have

\vspace{2mm}

(i). If $\sigma_*= \sigma_{lr}, \sigma_{SF0},\sigma_{su},
\sigma_{r}, \sigma_{a}, \sigma_{l}$ or $\sigma_{se}$, then
$\eta(\sigma_{*}(T))=\eta(\sigma(T))$.

\vspace{2mm}

(ii). If $\sigma_*=\sigma_{ab}, \sigma_{lb}, \sigma_{sb},
  \sigma_{rb}, \sigma_{aw}, \sigma_{lw}, \sigma_{sw},  \sigma_{rw}, \sigma_{w},
  \sigma_{SF+}, \sigma_{SF-}, \sigma_{le}, \sigma_{re}, \sigma_{e}, \sigma_{K}$ or
  $\sigma_{K_3}$, then $\eta(\sigma_{*}(T))=\eta(\sigma_{b}(T))$.

  \vspace{2mm}

(iii).
$\eta(\sigma_{D}(T))=\eta(\sigma_{des}(T))=\eta(\sigma_{rD}(T))=\eta(\sigma_{lD}(T))$.

\vspace{3mm}

\par \noindent{\bf Proof.} (i). Note that if $\sigma_*=\sigma_{lr}, \sigma_{SF0},
\sigma_{su}, \sigma_{r}, \sigma_{a}$ or $\sigma_{l}$, then
$\partial\sigma\subseteq
\sigma_{su}\cap\sigma_{a}\subseteq\sigma_{*}\subseteq\sigma$,
$\partial\sigma\subseteq \sigma_{se} \subseteq\sigma$ ([1]), so (i)
is proved.

\vspace{2mm}

\par (ii). If $\sigma_*=\sigma_{ab}, \sigma_{lb}, \sigma_{sb},
  \sigma_{rb},\sigma_{aw}, \sigma_{sw},  \sigma_{sw}, \sigma_{rw}, \sigma_{w},
  \sigma_{SF+}, \sigma_{SF-},
 \sigma_{le}, \sigma_{re}, \sigma_{e}, \sigma_{K}, \sigma_{K_3}$ or $\sigma_{SF+}$, it is well known that
$\partial\sigma_b\subseteq\sigma_{K}\subseteq \sigma_{*}
\subseteq\sigma_{b}$, so we have
$\eta(\sigma_{*}(T))=\eta(\sigma_b(T))$.

\vspace{2mm}

\par(iii). First, we prove that
$\eta(\sigma_{D}(T))=\eta(\sigma_{des}(T)).$ That $
\sigma_{des}(T)\subseteq\sigma_{D}(T)$ is clear. If $\lambda\in
\partial(\sigma_{D}(T))$, there exist $\{\lambda_n\}$ such that
$$\{\lambda_n\}\cap\sigma(T)=\emptyset\,\, \makebox{and}\,\,
\lambda_n\rightarrow \lambda.$$ If $\lambda\not\in\sigma_{des}(T),$
then $\sigma_{des}(T-\lambda)<\infty$ and hence by Lemma 2.10 that
$\lambda\not\in \sigma_{D}(T)$, which contradicts with $\lambda\in
\partial(\sigma_{D}(T))\subseteq\sigma_{D}(T)$. Thus it follows
that $\lambda\in\sigma_{des}(T)$ when $\lambda\in
\partial(\sigma_{D}(T))$,  so $
\partial(\sigma_{D}(T))\subseteq\sigma_{des}(T).$ Note that
$\sigma_{des}(T)\subseteq\sigma_{rD}(T)\subseteq\sigma_{D}(T)$,
thus,
 $\eta(\sigma_{D}(T))=\eta(\sigma_{rD}(T)).$ Moreover, since
 $\sigma_{rD}$ and $\sigma_{lD}$ are dual, we have
 $$\eta(\sigma_{lD}(T))=\eta(\sigma_{rD}(T^*))=\eta(\sigma_{D}(T^*))=\eta(\sigma_{D}(T)).$$
Therefore,
$\eta(\sigma_{D}(T))=\eta(\sigma_{des}(T))=\eta(\sigma_{rD}(T))=\eta(\sigma_{lD}(T))$.
 This proved (iii).

 \vspace{3mm}

The following theorems decide 18 kind spectra filling-in-hole
properties of $M_C$:

\vspace{3mm}

\vskip 0.2in\noindent {\bf Theorem 4.3.} If $\sigma_*=\sigma_{des},
\sigma_{su},  \sigma_{r},\sigma_{rb},\sigma_{sw}, \sigma_{rw}$ or
$\sigma_{re}$, then $\sigma_* $ has the generalized filling-in-hole
property.
  That is $$\sigma_{*}(A)\cup \sigma_{*}(B)=\sigma_{*}(M_C)  \cup W_{\sigma_{*}}(A,B,C),$$
where $W_{\sigma_{*}}(A,B,C)$ is contained in the union of some
holes in $\sigma_{*}(M_C)$.  In particular,

\vspace{2mm}

\par (i). $W_{\sigma_{des}}(A,B,C)\subseteq [(\sigma_{des}(A)\cap
    \sigma_{asc}(B))\setminus \sigma_{des}(B)]$ is  contained in the union of all holes in
    $\sigma_{des}(B).$

    \vspace{2mm}

\par(ii). If $\sigma_*=\sigma_{su}, \sigma_{r},
  \sigma_{rb}$ or $\sigma_{re}$, then  $W_{\sigma_{*}}(A,B,C)\subseteq [(\sigma_{*}(A)\cap
    \sigma_{*}(B^*))\setminus \sigma_{*}(B)]$ is  contained in  the union of all holes  in
    $\sigma_{*}(B).$

    \vspace{2mm}

\par(iii). If $\sigma_*=\sigma_{rw}$ or $\sigma_{sw}$, then  $W_{\sigma_{*}}(A,B,C)\subseteq [(\sigma_{*}(A)\cup
    \sigma_{*}(B^*))\setminus \sigma_{SF-}(B)]$ is  contained in  the union of all holes  in
    $\sigma_{*}(B).$

\medskip
\noindent {\em \bf Proof.} (i). It follows from Lemma 2.5 and the
proof of Lemma 2.6 (i) that
$$\sigma_{des}(A)\cup  \sigma_{des}(B)\supseteq\sigma_{des}(M_{C})
\supseteq \sigma_{des}(B).$$ This implies that $\sigma_{des}$ has
the generalized filling-in-holes property and
$$W_{\sigma_{des}}(A,B,C)\subseteq \sigma_{des}(A)\setminus
\sigma_{des}(B).$$  Moreover, note that Lemma 2.7, we can prove that
$$W_{\sigma_{des}}(A,B,C)\subseteq [(\sigma_{des}(A)\cap\sigma_{asc}(B))\setminus
\sigma_{\sigma_{des}}(B)].$$  To see this, let $\lambda\in
W_{\sigma_{des}}(A,B,C).$ Then $\sigma_{des}(A-\lambda)=\infty$ and
$\sigma_{des}(B-\lambda)<\infty.$ If
$\sigma_{asc}(B-\lambda)<\infty,$ then by Lemma 2.7 (ii) we know
that $\sigma_{des}(A-\lambda)<\infty$, which is a contradiction.
Thus $\lambda\in\sigma_{asc}(B)$. Next, we can claim that
$W_{\sigma_{des}}(A,B,C)$ is contained in the union of all holes in
$\sigma_{des}(B),$ that is, $W_{\sigma_{des}}(A,B,C)
\subseteq\eta(\sigma_{des}(B)).$ Otherwise, there exists $\lambda\in
W_{\sigma_{des}}(A,B,C) \setminus\eta(\sigma_{des}(B)).$ By
Proposition 4.2 we have that
$\eta(\sigma_{des}(B))=\eta(\sigma_D(B))$. Thus $\lambda\not\in
\eta(\sigma_{D}(B)).$ Furthermore, Lemma 2.7 tells us that
$des(A-\lambda)<\infty\Leftrightarrow des(M_C-\lambda)<\infty,$
which is a contradiction with the assumption that $\lambda\in
W_{\sigma_{des}}(A,B,C).$ Thus, it follows that
$W_{\sigma_{des}}(A,B,C)$ is contained in the union of all holes in
$\sigma_{des}(B).$

 \vskip 0.1in

 (ii). We only prove $\sigma_* =\sigma_{su}$ case. Note that $A$ and $B$ are surjective imply that $M_C$ is surjective, $M_C$ is surjective implies that $B$ is also surjective. So
 we have
$$\sigma_{su}(B)\subseteq\sigma_{su}(M_{C})\subseteq \sigma_{su}(A)\cup  \sigma_{su}(B)
.$$ This shows that $\sigma_{su}$ has the generalized
filling-in-hole property and
$$W_{\sigma_{su}}(A,B,C)\subseteq \sigma_{su}(A)\setminus
\sigma_{su}(B).$$ Next we claim that
$W_{\sigma_{su}}(A,B,C)\subseteq \sigma_{a}(B).$  If not, there
exists  $\lambda\in W_{\sigma_{su}}(A,B,C)\setminus\sigma_{a}(B).$
Combine this fact with the inclusion
$W_{\sigma_{su}}(A,B,C)\subseteq \sigma_{su}(A)\setminus
\sigma_{su}(B)$ proved above, we have that $B-\lambda$ is
invertible. By Lemma 2.7 it follows that $A-\lambda$ is surjective,
this is a contradiction. Thus $W_{\sigma_{su}}(A,B,C)\subseteq
\sigma_{a}(B)$, and so
$$W_{\sigma_{su}}(A,B,C)\subseteq
(\sigma_{su}(A)\cap\sigma_{a}(B))\setminus \sigma_{su}(B)=
(\sigma_{su}(A)\cap\sigma_{su}(B^*))\setminus \sigma_{su}(B).$$
 Similar to the proof of (1), we can get that $W_{\sigma_{su}}(A, B, C)$
is contained in  the union of all holes in $\sigma_{su}(B).$

 \vskip 0.1in

 (iii) is obvious by Proposition 4.2.

  \vspace{3mm}

 Duality, we have the following:

 \vskip 0.2in\noindent {\bf Theorem 4.4.} If $\sigma_*=\sigma_{l}, \sigma_{lb}, \sigma_{aw}, \sigma_{lw}$ or
 $\sigma_{le}$, then $\sigma_* $ has the generalized  filling-in-holes property. That is
$$\sigma_{*}(A)\cup \sigma_{*}(B)=\sigma_{*}(M_C)  \cup W_{\sigma_*}(A,B,C),$$
where $W_{\sigma_*}(A,B,C)$ is contained in the union of some holes
in $\sigma_{*}(M_C)$. In particular,

 \vspace{2mm}

\par(i). If $\sigma_*=\sigma_{l},\sigma_{lb}$ or $\sigma_{le}$, then  $W_{\sigma_*}(A,B,C)\subseteq [(\sigma_{*}(B)\cap
    \sigma_{*}(A^*))\setminus \sigma_{*}(A)]$ is  contained in  the union of all holes in
    $\sigma_{*}(A).$

     \vspace{2mm}

\par(ii). If $\sigma_*=\sigma_{aw}$ or $\sigma_{lw}$, then $W_{\sigma_*}(A,B,C)\subseteq [(\sigma_{*}(B)\cup
    \sigma_{*}(A^*))\setminus \sigma_{SF+}(A)]$ is contained in the union of all holes  in
    $\sigma_{*}(A).$

     \vspace{2mm}

 \vskip 0.2in\noindent {\bf Theorem 4.5.} If
$\sigma_*=\sigma_{lD},\sigma_{rD},\sigma_{lr},
\sigma_{K_3},\sigma_{K}$ or $\sigma_{SF_0}$, then $\sigma_*$ has the
convex filling-in-hole property. That is
$$\eta(\sigma_{*}(A)\cup \sigma_{*}(B))=\eta(\sigma_{*}(M_C)).$$
\medskip
\noindent {\em \bf Proof.} It follows from Lemma 3.1 of [12] that
$\eta(\sigma(A)\cup \sigma(B))=\eta(\sigma(M_C)),
\eta(\sigma_{b}(A)\cup \sigma_{b}(B))=\eta(\sigma_{b}(M_C)),$
$\eta(\sigma_{D}(A)\cup \sigma_{D}(B))=\eta(\sigma_{D}(M_C)).$
Combine those facts with Proposition 4.2 that it is easy to prove
the theorem.

\vspace{3mm}

We are also interesting in the following question:

 \vspace{3mm}

\medskip\par {\bf Open question 4.6.} Do other spectra of $M_C$ have the filling-in-hole properties ?

\large\section{Examples}

\par Now, we present examples to show that some conclusions about the spectra structure
or the filling-in-hole properties of $M_C$ are not true.

 \vspace{3mm}

The following example shows that for spectrum
$\sigma_{*}=\sigma_{a},\sigma_{l},\sigma_{SF+},\sigma_{le},\sigma_{aw}$,
$\sigma_{lw},\sigma_{ab}, \sigma_{lb}, \sigma_{su},$
$\sigma_{r},\sigma_{SF-},\sigma_{re}$,$\sigma_{sw},\sigma_{rw},\sigma_{sb}$
or $\sigma_{rb}$, it not only has not the equation (1) form, but
also has not the filling-in-hole property.

 \vspace{3mm}

\vskip 0.2in \noindent {\bf  Example 5.1} ([12]). Let $X$ be the
direct sum of countably many copies of $\ell^2:=\ell^2(N)$. Thus,
the elements of $X$ are the sequences $\{x_j\}_{j=1}^\infty$ with
$x_j \in \ell^2$ and $\sum_{j=1}^\infty \| x_j\|^2 < \infty$. Put
$Y=\ell^2$. Let $V$ be the forward shift on $\ell^2$:
$$V : \ell^2 \to \ell^2, \quad \{z_1,z_2, \ldots\} \mapsto \{0,z_1,z_2,\ldots\},$$
define the operators $A$ and $C$ by
$$ A: X \to X, \quad \{x_1, x_2, \ldots\} \mapsto \{Vx_1, Vx_2, \ldots\},$$
$$ C: Y\to X, \quad \{y_1, y_2, \ldots\} \mapsto \{y_1e_1, y_2e_1, \ldots\}.$$
where $e_1=\{1,0,0,\ldots\}$. Let $B=0$ and consider the operator
$$M_C=\left(\begin{array}{ll} A & C \\ 0 &B \end{array}\right) : X \oplus Y \to X \oplus
Y.$$

If
$\sigma_{*}=\sigma_{a},\sigma_{l},\sigma_{SF+},\sigma_{le},\sigma_{aw},\sigma_{lw},\sigma_{ab}$
or $\sigma_{lb}$, then

 \vspace{2mm}

\par (i). $\sigma_{*}(A)\cup\sigma_{*}(B)\cup(S(A^*)\cap S(B))
=\sigma_{*}(A)\cup\sigma_{*}(B)=\{\lambda\in
{\mathbb{C}}:\mid\lambda\mid=1\}\cup\{0\}.$

 \vspace{2mm}

\par (ii). $\sigma_{*}(M_C)\cup(S(A^*)\cap S(B))=\sigma_{*}(M_C)=\{\lambda\in {\mathbb{C}}:\mid\lambda\mid=1\}.$

 \vspace{2mm}

\par (iii). $(\sigma_{*}(A)\cup\sigma_{*}(B))\setminus\sigma_{*}(M_C)=\{0\}.$

 \vspace{2mm}

Thus, it follows from (i) and (ii) that equation (1) does not hold
for spectrum
$\sigma_{*}=\sigma_{a},\sigma_{l},\sigma_{SF+},\sigma_{le},\sigma_{aw}$,
$\sigma_{lw},\sigma_{ab}$ or $\sigma_{lb}$. By duality, we can also
show that equation (1) does not hold for spectrum
$\sigma_{*}=\sigma_{su},\sigma_{r},\sigma_{SF-},\sigma_{re}$,
$\sigma_{sw},\sigma_{rw},\sigma_{sb}$ or $\sigma_{rb}$. Moreover,
from (iii) we knew that none of the above 16 kind spectra has the
filling-in-hole property.

 \vspace{3mm}

This following example shows that ascent spectrum $\sigma_{asc}$ has
not equations (1) to (9) form.

\vskip 0.2in\noindent {\bf  Example 5.2.} Let $X=Y=\ell^2$ and
$\{e_n\}_{n\geq 1}$ be a basis of $\ell^2$. Define $A, B, C \in
B(\ell^2)$ by

$$A e_{i}=\frac{1}{i}e_{2i}\,\,\, \makebox{for} \,\, i=1,2,\cdots,$$
$$B e_1= 0, B e_i= \frac{1}{i}e_{i-1} \,\,\,\,\,\makebox{for} \,\,i=2,3,\cdots,$$
$$C e_{i}= e_{2i-1} \,\,\,\,\makebox{for} \,\,i=1,2,3,\cdots.$$
Then we have
\begin{eqnarray*}
& \{0\}=\sigma_{asc}(A)\cup \sigma_{asc}(B) &=\sigma_{asc}(A)\cup
\sigma_{asc}(B) \cup S(A)\cup S(A^*)\cup S(B)\cup S(B^*)\\
&&\not= \sigma_{asc}(M_{C})\cup S(A)\cup S(A^*)\cup S(B)\cup S(B^*)\\
& &= \sigma_{asc}(M_{C})=\emptyset.\\
\end{eqnarray*}

 \vspace{3mm}

In [16], the authors claimed that
$$(\sigma_{ab}(A)\cup   \sigma_{ab}(B))\setminus
\sigma_{ab}(M_C )  \subseteq  S(A^*)\cap \sigma_D(B),$$where
$\sigma_D(B)$ was denoted by $F(B)$, which implies that
$$(\sigma_{ab}(A)\cup \sigma_{ab}(B))\cup (S(A^*)\cap
\sigma_D(B))=\sigma_{ab}(M_C ) \cup  (S(A^*)\cap
\sigma_D(B))\eqno(13).$$

\vspace{3mm}

The following example shows that neither the claim nor equation (13)
is true.

\vskip 0.2in\noindent {\bf Example 5.3.} Let $X=Y=\ell^2$ and
$$A   \{x_1, x_2, \ldots\} \mapsto \{ x_1, 0,  x_2, 0,\ldots\},$$
$$B  \{x_1, x_2, \ldots\} \mapsto \{0, 0, 0, \ldots\},$$
$$C  \{x_1, x_2, \ldots\} \mapsto \{0, x_1, 0, x_2,
\ldots\}.$$ Then
$$\sigma_{ab}(A) =\{\lambda\in {\mathbb{C}}:\mid\lambda\mid=1\},   \sigma_{ab}(B) =\{0\}, \sigma_{ab}(M_C)
 =\{\lambda\in {\mathbb{C}}:\mid\lambda\mid=1\},$$ $$S(A^*)\cap \sigma_D(B)=\emptyset, \sigma_{ab}(M_C)\not=\sigma_{ab}(A)\cup
  \sigma_{ab}(B).$$ So $$(\sigma_{ab}(A)\cup   \sigma_{ab}(B))\setminus \sigma_{ab}(M_C )
\not \subseteq  S(A^*)\cap  \sigma_D(B).$$

 \vspace{3mm}

The following example shows that for spectrum
$\sigma_{*}=\sigma_{aw},\sigma_{lw},\sigma_{sw}$ or $\sigma_{rw}$,
it has not equations (4) and (5) form.

\vskip 0.2in \noindent {\bf  Example 5.4.}  Let $X=Y=\ell^2$ and
define $T, S, C \in B(\ell^2)$ by
 $$T \{x_1, x_2, x_3, \ldots\}\mapsto \{0, x_1, x_2,\ldots\},$$
$$S  \{x_1, x_2, x_3,  \ldots\} \mapsto \{x_2, x_4, x_6, \ldots,\},$$
$$C \{x_1, x_2, x_3, \ldots\} \mapsto \{0, 0, 0,  \ldots\}.$$

(i). Put $A=S, B=T^2$. Then $A^*$ and $B$ have the SVEP and
$$\sigma_{aw}(A)\cup\sigma_{aw}(B)\cup S(A^*)\cup S(B)=\sigma_{lw}(A)\cup \sigma_{lw}(B)\cup S(A^*)\cup
S(B)=\{\lambda\in {\mathbb{C}}:\mid\lambda\mid\leq1\},$$ and $$
\{\lambda\in {\mathbb{C}}:\mid\lambda\mid=1\} =\sigma_{aw}(M_C)\cup
S(A^*)\cup S(B) =\sigma_{lw}(M_C)\cup S(A^*)\cup S(B).$$ So, when
$\sigma_{*}=\sigma_{aw}$ or $\sigma_{lw}$, we have
$$\sigma_{*}(A)\cup\sigma_{*}(B)=\sigma_{*}(A)\cup\sigma_{*}(B)\cup
S(A^*)\cup S(B)\not=\sigma_{*}(M_C)\cup S(A^*)\cup
S(B)=\sigma_{*}(M_C).$$

 \vspace{3mm}

(ii). Put $A=S^2, B=T$. Then $A^*$ and $B$ have the SVEP and
$$\sigma_{sw}(A)\cup\sigma_{sw}(B)\cup S(A^*)\cup
S(B)=\sigma_{rw}(A)\cup \sigma_{rw}(B)\cup S(A^*)\cup
S(B)=\{\lambda\in {\mathbb{C}}:\mid\lambda\mid\leq1\}$$ and
$$\{\lambda\in {\mathbb{C}}:\mid\lambda\mid=1\}
=\sigma_{sw}(M_C)\cup S(A^*)\cup S(B) =\sigma_{rw}(M_C)\cup
S(A^*)\cup S(B).$$ Thus, when $\sigma_{*}=\sigma_{sw}$ or
$\sigma_{rw}$, we have
$$\sigma_{*}(A)\cup\sigma_{*}(B)=\sigma_{*}(A)\cup\sigma_{*}(B)\cup
S(A^*)\cup S(B)\not=\sigma_{*}(M_C)\cup S(A^*)\cup
S(B)=\sigma_{*}(M_C).$$

 \vspace{3mm}

The following example shows that for spectrum
$\sigma_{*}=\sigma_{K}, \sigma_{K_3}, \sigma_{SF0}$ or
$\sigma_{lr}$, it not only has not equations (4) and (5) form, but
also has not generalized filling-in-hole property.

\vspace{3mm}

\vskip 0.2in \noindent {\bf  Example 5.5.}  Let $X=Y=\ell^2$  and
define $A, B, C \in B(\ell^2)$ by
$$A \{x_1, x_2, x_3, \ldots\} \mapsto \{ x_2, x_4, x_6,\ldots\},$$
$$B \{x_1, x_2, x_3,  \ldots\} \mapsto \{0, x_1, 0, x_2, \ldots,\},$$
$$C \{x_1, x_2, x_3, \ldots\} \mapsto \{0, 0, 0,  \ldots\}.$$
It is easy to show that $A^*$ and $B $ have the SVEP and if
$\sigma_{*}=\sigma_{K}, \sigma_{K_3}, \sigma_{SF0}$ or
$\sigma_{lr}$, then $\sigma_{*}(A)\cup \sigma_{*}(B) =\{\lambda\in
{\mathbb{C}}:\mid\lambda\mid=1\}$, $ \sigma_{*}(M_C)
 =\{\lambda\in {\mathbb{C}}:\mid\lambda\mid\leq1\}.$ Thus, when $\sigma_{*}=\sigma_{K}, \sigma_{K_3}, \sigma_{SF0}$ or $\sigma_{lr}$, we have $$\sigma_{*}(A)\cup \sigma_{*}(B)\cup S(A^*)\cup
S(B)\not = \sigma_{*}(M_C)\cup S(A^*)\cup S(B)$$ and
$$\sigma_{*}(A)\cup \sigma_{*}(B)\not\supseteq
\sigma_{*}(M_C).$$

 \vspace{3mm}

The following example shows that the inclusions in Theorem 4.1 may
be strict in general.

 \vspace{3mm}

\vskip 0.2in \noindent {\bf  Example 5.6.}
 Let $X_n$ be a complex $n$ dimensional
Hilbert space. Define $T, S, C_3\in B( \ell^2)$ by
$$T \{x_1, x_2,x_3,\cdots\}=\{0, x_1, x_2, x_3,\cdots\},$$
$$S \{x_1, x_2,x_3,\cdots\}=\{x_2, x_3, x_4,\cdots\},$$
$$C_3: =I-TS.$$

(i). Put $A=T$,
  $C=\left(
       \begin{array}{cc}
         C_3&0
       \end{array}
     \right)$   $:  \ell^2 \oplus X_n \longrightarrow  \ell^2$,
$B=\left(
       \begin{array}{cc}
         S&0\\
         0&0
       \end{array}
     \right)$   $:  \ell^2 \oplus X_n\longrightarrow  \ell^2 \oplus
     X_n$. Then $W_{\sigma_{}}(A,B,C)=\{\lambda\in {\mathbb{C}}:0<\mid \lambda
\mid < 1$ $ \},$ $W_{\sigma_{b}}(A,B,C)=\{\lambda\in
{\mathbb{C}}:\mid \lambda \mid < 1$ $ \}$, thus
$W_{\sigma_{}}(A,B,C) \not=W_{\sigma_{b}}(A,B,C)$.

(ii). Put $A=T$, $C=\left(
       \begin{array}{cc}
         C_3&0
       \end{array}
     \right)$   $:  \ell^2 \oplus  \ell^2 \longrightarrow  \ell^2$,
$B=\left(
       \begin{array}{cc}
         S&0\\
         0&0
       \end{array}
     \right)$   $:  \ell^2 \oplus  \ell^2\longrightarrow  \ell^2 \oplus
     \ell^2$. Then  $W_{\sigma_{b}}(A,B,C)=\{\lambda\in {\mathbb{C}}:0<\mid \lambda
\mid < 1$ $ \},$ $W_{\sigma_{D}}(A,B,C)=\{\lambda\in
{\mathbb{C}}:\mid \lambda \mid < 1$ $ \}$, thus
$W_{\sigma_{b}}(A,B,C) \not=W_{\sigma_{D}}(A,B,C)$.

(iii). Put $A=T$,  $C=0$, $B=S$. Then
$W_{\sigma_{b}}(A,B,C)=\emptyset,$
$W_{\sigma_{w}}(A,B,C)=\{\lambda\in {\mathbb{C}}:\mid \lambda \mid <
1\}$, thus $W_{\sigma_{b}}(A,B,C) \not=W_{\sigma_{w}}(A,B,C)$.

\vspace{3mm}
 The following example shows that for $W_{\sigma_{e}}(A,B,C)$, $W_{\sigma_{w}}(A,B,C)$ and $W_{\sigma_{D}}(A,B,C)$,
 there are no inclusion relationship among them.

  \vspace{3mm}

 \vskip 0.2in \noindent {\bf  Example 5.7.} Let
  $T, S, T_1, S_1, C_1, C_2, C_3\in B( \ell^2)$ be defined by
 $$T \{x_1,x_2,x_3,\cdots\}=\{0,x_1,x_2,x_3,\cdots\},$$
 $$S \{x_1,x_2,x_3,\cdots\}=\{x_2,x_3,x_4,\cdots\},$$
 $$T_1 \{x_1,x_2,x_3,\cdots\}=\{0,x_1,0,x_2,0, x_3,\cdots\},$$
 $$S_1 \{x_1,x_2,x_3,\cdots\}=\{x_2,x_4,x_6,\cdots\},$$
 $$C_1 \{x_1,x_2,x_3,\cdots\}=\{x_1,0,x_3,0,x_5,\cdots\},$$
 $$C_2 \{x_1,x_2,x_3,\cdots\}=\{0,0,x_1,0,x_3,0,x_5,\cdots\},$$
 $$C_3 =I-TS.$$

(i).  If $A=T, B=S, C=C_3$, then $W_{\sigma_{w}}(A,B,C) \not
\subseteq W_{\sigma_{e}}(A,B,C).$

(ii). If $A=T, B=S, C=0$, then $W_{\sigma_{w}}(A,B,{C}) \not
\subseteq W_{\sigma_D}(A,B,{C}).$

(iii). Let $A=T_1,$ $B= \left(
       \begin{array}{cc}
        S_1&0\\
         0&S\\
       \end{array}
            \right)$$:  \ell^2 \oplus  \ell^2 \longrightarrow   \ell^2 \oplus  \ell^2 $
   and   $C= \left(
       \begin{array}{cc}
        C_1&0
       \end{array}
            \right)$$:   \ell^2 \oplus  \ell^2  \longrightarrow   \ell^2. $
 Then $W_{\sigma_{e}}(A,B,{C}) \not \subseteq W_{\sigma_{D}}(A,B,{C}).$

(iv). If $A=T_1, B=S_1, C=C_2$, then $W_{\sigma_{e}}(A,B,{C}) \not
\subseteq W_{\sigma_{w}}(A,B,{C}).$

(v). If $A=T$, $ B = \left(
       \begin{array}{cc}
        S&0\\
         0&0\\
       \end{array}
            \right)$$:  \ell^2 \oplus  \ell^2 \longrightarrow   \ell^2 \oplus  \ell^2 $
   and   $C  = \left(
       \begin{array}{cc}
        C_3&0
       \end{array}
            \right)$$:   \ell^2 \oplus  \ell^2  \longrightarrow   \ell^2,
            $ then $W_{\sigma_{D}}(A,B,{C}) \not \subseteq W_{\sigma_{e}}(A,B,{C}), W_{\sigma_{D}}(A,B,{C}) \not \subseteq W_{\sigma_{w}}(A,B,{C}).$

             \vspace{5mm}

Note that $\sigma_{}(T)\supseteq\sigma_{b}(T)
\supseteq\sigma_{D}(T)\supseteq acc\sigma(T)$ is well known, so we
have $\eta(\sigma_{}(T))\supseteq\eta(\sigma_{b}(T))
\supseteq\eta(\sigma_{D}(T))\supseteq\eta(acc\sigma(T)).$

 \vspace{3mm}

 The following example shows that the above inclusions may be strict.

\vskip 0.2in \noindent {\bf  Example 5.8.} Let $X_n$ be a $n$
dimensional complex Hilbert space. Define operators $A \in B(
\ell^2)$ by $$A\{x_1, x_2, x_3, \cdots\} = \{0,\frac {1}{2}x_1,
\frac {1}{3}x_2, \frac {1}{4}x_3, \cdots\}.$$ Then
$\sigma(A)=\sigma_{D}(A)=\sigma_{des}(A)=\{0\}$ and
$acc\sigma(A)=\emptyset$. If consider operator $$T=\left(
       \begin{array}{ccc}
        0&  0 &0\\
         0&3I&0\\
         0 &0&5+A
       \end{array}
     \right): X_n\oplus  \ell^2\oplus  \ell^2  \longrightarrow X_n\oplus  \ell^2\oplus  \ell^2,$$ we have $\eta(\sigma_{}(T))=\sigma_{}(T) =\{0,3,5\}$,
     $\eta(\sigma_{b}(T))=\sigma_{b}(T) =\{3,5\}$, $\eta(\sigma_{D}(T))=\sigma_{D}(T) =\{5\}$, $\eta(acc\sigma(T))=acc\sigma(T)
     =\emptyset$. Thus, $$\eta(\sigma_{}(T))\not=\eta(\sigma_{b}(T))\not=\eta(\sigma_{D}(T))\not=
  \eta(acc\sigma(T)).$$

   \vspace{3mm}

The following example shows that spectra $\sigma_{rD}$ and
$\sigma_{lD}$ have not the generalized filling-in-hole property.

 \vspace{3mm}

\vskip 0.2in \noindent {\bf  Example 5.9.} Let $X=Y=\ell^2$. Define
$A=0$ and $B, C$ by
$$B\{x_1,x_2,x_3,\cdots\}=\{x_2,x_4,x_6,\cdots\},$$
$$C\{x_1,x_2,x_3,\cdots\}=\{x_1, \frac {1}{\sqrt2}x_3,\frac
{1}{\sqrt3}x_5, \cdots\}.$$ Then $0\not \in \sigma_{rD}(A)\cup
\sigma_{rD}(B)$ but $0\in \sigma_{rD}(M_C)$. So
$\sigma_{rD}(M_C)\not\subseteq\sigma_{rD}(A)\cup\sigma_{rD}(B).$
Since $\sigma_{rD}$ and $\sigma_{lD}$ are dual (see [26]), thus,
neither $\sigma_{rD}$ nor $\sigma_{lD}$ has the generalized
filling-in-holes property.

 \vspace{5mm}

Let $H, K$ be Hilbert spaces, $(A, B)\in B(H)\times B(K), C\in B(K,
H)$. If $A\in B(H)$, let $A^{*}$ denote the adjoint operator of $A$
and $\sigma_p (A)$ denote the point spectrum of $A$. In [3], the
authors claimed that
$$\eta(\sigma_{se}(A)\cup \sigma_{se}(B))= \eta(
\sigma_{se}(M_C)),$$ More precisely,
$$\sigma_{se}(A)\cup  \sigma_{se}(B) \cup
(\overline{\sigma_{p}(A^{*})}\cap
\sigma_{p}(B))=\sigma_{se}(M_C)\cup
      W,\eqno(14)$$
where $W$ is the union of some holes in $\sigma_{se}(M_C)$ which
happen to be subsets of $\overline{\sigma_{p}(A^{*})}\cup
\sigma_{p}(B)$. The following example shows that equation (14) is
not true.

\vspace{3mm}

\vskip 0.2in \noindent {\bf  Example 5.10.} Let  $X, Y, A, C$  be
defined as in Example 5.1. and $B\in B(Y)$ be defined by
$$ B: \{y_1, y_2, \ldots\} \mapsto \{0,\frac {1}{2}y_1, \frac
{1}{3}y_2,\frac {1}{4}y_3, \cdots\}.$$ Consider the operator
$$M_C=\left(\begin{array}{ll} A & C \\ 0 & B \end{array}\right) : X \oplus Y \to X \oplus Y.$$

Then we have

 \vspace{2mm}

(i).
$\sigma_{se}(M_C)=\sigma_{se}(A)=\{\lambda:\mid\lambda\mid=1\}$,

 \vspace{2mm}

(ii). $\sigma_{}(B)=\sigma_{{se}}(B)=\{0\},$
$\sigma_{p}(B)=\emptyset,$

 \vspace{2mm}

(iii). $\overline{\sigma_p( A^*)}\cap \sigma_p (B)=\emptyset.$

 \vspace{3mm}

Thus $
W=(\sigma_{{se}}(A)\cup\sigma_{{se}}(B)\cup(\overline{\sigma_p
(A^*)}\cap \sigma_p (B))\setminus\sigma_{{se}}(M_C)=\{0\},$ so $W$
is just a point but not an open set. This showed that the above
conclusion is not true.

\medskip\noindent {{ \large \bf Acknowledgment}}
\par  The authors are grateful to Doctor Qiaofen Jiang for  the valuable
suggestions on Lemma 2.7 and  Theorem 3.1.


\begin{thebibliography}{S2}



\bibitem{1} P. Aiena. Fredholm and Local Spectral Theory, with Applications to Multipliers, Kluwer Academic Publishers, 2004.
\bibitem{2} M. Barraa and M. Boumazgour. A note on the  spectrum of an upper triangular operator matrix,  Proc. Amer. Math. Soc., 131(2003), 3083-3088.
\bibitem{3} M. Barraa, M. Boumazgour. On the perturbations of spectra of upper triangular operator matrices,  J. Math. Anal. Appl., 347(2008), 315-322.
\bibitem{4} M. Berkani. Index of B-Fredholm operators and generalization of a Weyl theorem, Proc. Amer. Math. Soc., {130}(2002), 1717-1723.
\bibitem{5} M. Berkani. On a class of quasi-Fredholm operators, Integr. Equ. Oper. Theory, 34(1999), 244-249.
\bibitem{6} C. Benhida, E. H. Zerouali and H. Zguitti. Spectral properties of upper-triangular block operators, Acta Sci. Math. (Szeged), 71(2005), 607-616.
\bibitem{7} M. Burgos, A. Kaidi, M. Mbekhta and  M. Oudghiri. The descent spectrum and perturbations, J. operator Theory, 56(2006), 259-271.
\bibitem{8} X. H. Cao. Browder spectra for upper triangular operator   matrices,\ J.  Math.  Anal. Appl., 342(2008), 477-484.
\bibitem{9} X. H. Cao, M.Z. Guo, B. Meng. Drazin  spectrum and Wely's theorem for operator matrices, J. Math.  Res. Exposition,  26( 2006), 413-422.
\bibitem{10} X. H. Cao, M.Z. Guo, B. Meng. Semi-Fredholm spectrum and Weyl's theory for operator matrices,\ Acta Math. Sinica, 22(2006), 169-178.
\bibitem{11} X. H. Cao, M.Z. Guo, B. Meng. Weyl's theorem for upper triangular operator matrices,  Linear Algebra  Appl., 402(2005), 61-73.
\bibitem{12} X. L. Chen, S. F. Zhang, H. J. Zhong. On the filling in holes problem of  operator matrices,  Linear Algebra  Appl., (2008), doi:10.1016/j.laa.2008.08.022.
\bibitem{13} D. S. Djordjevi\'{c}. Perturbations  of spectra  of operator matrices, J. Operator Theory, 48(2002), 467-486.
\bibitem{14} S. V. Djordjevi\'{c}, Y. M. Han. A note on Weyl's theorem for operator matrices, Proc. Amer. Math. Soc., 131(2002), 2543-2547.
\bibitem{15} S. V. Djordjevi\'{c}, Y. M. Han.  spectral continuity for operator matrices,\ Glasg.\  Math. J., 43(2001), 487-490.
\bibitem{16} S. V. Djordjevi\'{c},  H. Zguitti. Essential point   spectra of operator matrices though local spectral theory, J.  Math.  Anal.  Appl., 338(2008), 285-291.
\bibitem{17} H. K. Du, J. Pan. Perturbation of  spectrums  of $2\times2$ operator matrices,  Proc. Amer. Math. Soc., 121(1994), 761-766.
\bibitem{18} H. Elbjaoui, E. H. Zerouali. Local spectral theory for $2\times2$ operator matrices, Int. J. Math. Math. Sci., 42(2003), 2667-2672.
\bibitem{19} J. K. Han, H. Y. Lee, W. Y. Lee. Invertible completions of $2\times2$ upper triangular operator matrices.\  Proc. Amer. Math. Soc., 128(1999), 119-123.
\bibitem{20} I. S. Hwang, W. Y.  Lee. The boundedness  below of $2\times2$ upper triangular operator matrices, Integr. Equ. Oper. Theory, 39(2001), 267-276.
\bibitem{21} W. Y. Lee. Weyl's theorem for operator matrices, Integr. equ. oper. theory, 32(1998), 319-331.
\bibitem{22} W. Y. Lee. Weyl spectra of operator matrices,  Proc. Amer. Math. Soc., 129(2000), 131-138.
\bibitem{23} Y. Li, X. H. Sun, H. K. Du. A note on the left essential spectra of operator matrices, Acta Math. Sinica, 23(2007), 2235-2240.
\bibitem{24} Y. Li, X. H. Sun, H. K. Du. The intersection  of  left(right) spectra  of 2 $\times$ 2 upper triangular operator  matrices, Linear Algebra Appl., 418(2006), 112-121.
\bibitem{25} Y. Li, H. K. Du. The intersection  of essential approximate point spectra  of operator matrices, J. Math.  Anal.  Appl., 323(2006), 1171-1183.
\bibitem{26} M. Merkhta, V. M\"{u}ller. On the axiomatic theory of spectrum II., Studia Math., 119(1996), 129-147.
\bibitem{27} C. Schmoeger. Perturbation properties of some classes of operators,  Rendiconti di Matematica,  Serie VII  Volume 14 Roma, (1994), 533-541.
\bibitem{28} E. H. Zerouali, H. Zguitti. Perturbation of spectra of operator matrices and local spectral theory, J.  Math.  Anal.   Appl., 324(2006), 992-1005.
\bibitem{29} H. Y. Zhang, H. K. Du. Browder spectra of upper-triangular operator matrices,  J.  Math.  Anal. Appl., 323(2006), 700-707.
\bibitem{30} S. F. Zhang, H. J. Zhong. A note of Browder spectrum of operator matrices, J. Math. Anal. Appl., 344(2008), 927-931 .
\bibitem{31} S. F. Zhang, H. J. Zhong, Q. F. Jiang. Drazin  spectrum of operator matrices on the Banach space, Linear Algebra  Appl., 429(2008), 2067-2075.
\bibitem{32} Y. N. Zhang, H. J. Zhong, L. Q. Lin. Browder spectra and essential spectra of operator matrices, Acta Math. Sinica, 24 (2008), 947-954.

\end{thebibliography}
\end{document}